\documentclass[twocolumn]{autart}    

\usepackage{amsmath,amssymb,amsfonts}
\usepackage{algorithm,algorithmic}
\usepackage{graphicx}
\usepackage{textcomp}

\usepackage{color}

\usepackage{lipsum}
\usepackage{epstopdf}

\usepackage{amsopn}

\usepackage{enumerate}
\usepackage{cancel}
\usepackage{mathrsfs}
\usepackage{mathdots}
\usepackage{euscript}
\usepackage{amscd}
\usepackage{placeins}
\usepackage{natbib}

\newtheorem{theorem}{Theorem}

\newtheorem{lemma}{Lemma}
\newtheorem{proposition}{Proposition}
\newtheorem{assumption}{Assumption}

\theoremstyle{definition}

\theoremstyle{remark}
\newtheorem{remark}{Remark}

\newcommand{\bmat}{\left[ \begin{matrix}}
	\newcommand{\emat}{\end{matrix} \right]}
\newcommand{\innerprod}[2]{\langle{#1},\,{#2}\rangle}

\DeclareMathOperator{\trace}{tr}
\DeclareMathOperator{\argmax}{argmax}

\DeclareMathOperator{\E}{{\mathbb E}}
\newcommand{\Rbb}{\mathbb R}

\newcommand{\Jbb}{\mathbb J}
\newcommand{\Cbb}{\mathbb C}

\newcommand{\Zbb}{\mathbb Z}

\newcommand{\Tbb}{\mathbb T}

\newcommand{\yb}{\mathbf  y}
\newcommand{\sbf}{\mathbf  s}  
\newcommand{\zb}{\mathbf  z}

\newcommand{\fb}{\mathbf  f}
\newcommand{\gb}{\mathbf  g}

\newcommand{\nb}{\mathbf  n}

\newcommand{\tb}{\mathbf t}
\newcommand{\kb}{\mathbf k}
\newcommand{\lb}{\boldsymbol{\ell}}
\newcommand{\oneb}{\mathbf 1}
\newcommand{\zerob}{\mathbf 0}

\newcommand{\Fb}{\mathbf F}

\newcommand{\Nb}{\mathbf N}

\newcommand{\Qb}{\mathbf Q}

\newcommand{\Yb}{\mathbf Y}

\newcommand{\thetab}{\boldsymbol{\theta}}
\newcommand{\omegab}{\boldsymbol{\omega}}
\newcommand{\alphab}{\boldsymbol{\alpha}}

\newcommand{\Sigmab}{\boldsymbol{\Sigma}}

\newcommand{\zetab}{\boldsymbol{\zeta}}

\newcommand{\Hfrak}{\mathfrak{H}}

\newcommand{\Sfrak}{\mathfrak{S}}

\newcommand{\Lscr}{\mathscr{L}}

\newcommand{\Dscr}{\mathscr{D}}
\newcommand{\Mscr}{\mathscr{M}}

\newcommand{\Lcal}{\mathcal{L}}

\renewcommand{\d}{\mathrm{d}}
\newcommand{\p}{\mathrm{p}}
\newcommand{\F}{\mathrm{F}}

\begin{document}

\begin{frontmatter}

\title{M$^2$-Spectral Estimation:\\ A Relative Entropy Approach\thanksref{footnoteinfo}} 

\thanks[footnoteinfo]{This work was supported by the SID project ``A Multidimensional and Multivariate Moment Problem Theory for Target Parameter Estimation in Automotive Radars'' (ZORZ\_SID19\_01) funded by the Department of Information Engineering of the  University of Padova. The first author was also partially supported by the ``Hundred-Talent Program'' of the Sun Yat-sen University. Corresponding author B.~Zhu. Tel. +86 14748797525. 
Fax +86(20) 39336557.}

\author[SYSU]{Bin Zhu}\ead{zhub26@mail.sysu.edu.cn},    
\author[Unipd]{Augusto Ferrante}\ead{augusto@dei.unipd.it},               
\author[KTH]{Johan Karlsson}\ead{johan.karlsson@math.kth.se},  
\author[Unipd]{Mattia Zorzi}\ead{zorzimat@dei.unipd.it} 

\address[SYSU]{School of Intelligent Systems Engineering, Sun Yat-sen University, Waihuan East Road 132, 510006 Guangzhou, China}
\address[Unipd]{Department of Information Engineering, University of Padova, Via Gradenigo 6/B, 35131 Padova, Italy}  
\address[KTH]{Division of Optimization and Systems Theory, Department of Mathematics, KTH Royal Institute of Technology, 10044 Stockholm, Sweden}             


\begin{abstract}                          
This paper deals with M$^2$-signals, namely {\em multivariate} (or vector-valued) signals defined over a {\em multidimensional} domain. In particular, we propose an optimization technique to solve the covariance extension problem for stationary random vector fields. The multidimensional Itakura-Saito distance is employed as an optimization criterion to select the solution among the spectra satisfying a finite number of moment constraints. In order to avoid technicalities that may happen on the boundary of the feasible set, we deal with the discrete version of the problem where the multidimensional integrals are approximated by Riemann sums. The spectrum solution is also discrete, which occurs naturally when the underlying random field is periodic. We show that a solution to the discrete problem exists, is unique and depends smoothly on the problem data. Therefore, we have a well-posed problem whose solution can be tuned in a smooth manner. Finally, we have applied our theory to the target parameter estimation problem in an integrated system of automotive modules.  Simulation results show that our spectral estimator has promising performance.
\end{abstract}

\end{frontmatter}

\section{Introduction}\label{sec:intro}


Moment problems are ubiquitous in the areas of systems and control, where the moment conditions dictate system properties that needs to be satisfied. In this paper we consider a spectral estimation problem where the moment conditions ensure that the system covariances coincide with measured values \citep{stoica2005spectral,LP15}. The latter is known as rational covariance extension problem, which was initially proposed in \citet{Kalman}: 
For a given (scalar) partial covariance sequence 
of $n$ elements, 
determine all infinite  extensions such that the corresponding spectral density is nonnegative and rational with degree bounded by $n$. 
Note that the bound on the degree is a non-convex constraint, but it  naturally gives a bound on the complexity of the resulting system.

The solution to this problem was presented in \citet{byrnes1995acomplete} (see also \citet{georgiou_thesis} for an existence result) and  led to a convex optimization approach \citep{byrnes1998aconvex}
where the extension
is the maximizer of an entropy functional. This way of selecting extensions as maximizers of suitable functionals has been extensively studied, in particular in the unidimensional and univariate setting \citep{georgiou1999theinterpolation,byrnes2000anewapproach,byrnes2001cepstral,byrnes2002identifyability,enqvist2004aconvex}. Then, these spectral estimation paradigms have been extended to other types of functionals, which also typically come with guaranteed upper bounds of the degree of the extensions \citep{enqvist2008minimal,ferrante2008hellinger,zorzi2014anewfamily,Z-14rat}. More precisely, the power spectral density matches the partial covariance sequence and minimizes a pseudo-distance with respect to a prior spectral density which represents the a priori information on the system.
Several matrix-valued versions have also been considered \citep{georgiou2006relative,blomqvist2003matrix,ramponi2009aglobally,Z-15,FMP-12,zorzi2015interpretation,Pavon-Ferrante-SIAM-R-13,Zhu-Baggio-19,zhu2020well}.  
It is worth noting that the Nevanlinna-Pick interpolation problem is a special case of this framework, and this fact has been useful for applying the theory to control design \citep{nagamune2005sensitivity,takyar2008weight,karlsson2010theinverse,kergus2019reference}.

Most of these works are on dynamical systems in one variable (i.e. unidimensional systems), typically representing time. However, many problems in systems and control are inherently multidimensional \citep{bose2003multidimensional}. Multidimensional systems theory has been applied to many different problems, for example random Markov fields \citep{levy1990modeling}, image processing \citep{ekstrom1984digital} and target parameter estimation in radar applications \citep{rohling2012continuous,engels2014target}. Interest has therefore also been directed towards  multidimensional versions of the rational covariance extension problem \citep{georgiou2006relative,georgiou2005solution,ringh2015multidimensional,KLR-16multidimensional}. Most of the aforementioned works deal with the multidimensional and univariate case. However, there are situations in which the model is multidimensional and multivariate, say M$^2$. An example  
is given by the integrated system of automotive modules 
proposed in \citet{ZFKZ2019fusion}: the latter is composed by a certain number of uniform linear arrays (ULAs) of receive antennas sharing one common transmitter.

A natural approximation of the rational covariance extension problem is to restrict the support of the function to  a discrete grid. 
This was studied in  \citet{LPcirculant-13} for the unidimensional and univariate case, and is also called the  circulant rational covariance extension problem since it can be viewed as  limiting the  stationary process to be periodic with period $N$. A unidimensional and multivariate extension has been considered in \citet{lindquist2013onthemultivariate}, while the multidimensional and univariate case has been addressed in \citet{ringh2015multidimensional}. However the multidimensional and multivariable case has not been completely addressed yet.

In this paper, we consider the multidimensional and multivariable (M$^2$) version of the circulant rational covariance extension problem. This discrete version of the problem allows to avoid technicalities that may happen on the boundary of the feasible set. A natural choice of functional is the Itakura-Saito divergence \citep{enqvist2008minimal}, since it can be extended to matrix valued spectra and also allows for incorporating the a priori information \citep{FMP-12}. Moreover, as argued in \citet{FMP-12} in the $1$-d case, the IS-distance leads to a solution with low complexity. More precisely, the linear filter determined by the resulting spectrum has an a priori bounded McMillan degree, and the bound is as good as the one in the scalar case \citep{byrnes1995acomplete,byrnes1997partial}.
Thus, this leads to multivariate and multidimensional spectral analysis where information 
in terms of the covariances and the prior spectrum is fused in order to improve the estimates of the spectrum. Finally, we utilize the theory for parameter estimation in an integrated system of two automotive modules, and the numerical examples suggest that our method gives higher accuracy and robustness compared to the traditional periodogram-based method which represents the most straightforward way to compute an estimator of the spectrum from the data.

The outline of the paper is as follows. In Section \ref{sec:problem} we formulate the optimization problem.  In Section \ref{sec:dual} we prove the existence and uniqueness of the solution to the problem by means of duality theory. In Section \ref{sec:well-posed} we show that the solution depends continuously on the problem data. Then, we introduce the corresponding M$^2$ spectral estimator in Section \ref{sec:spec_est}, where we also provide a method to compute the covariance lags which guarantee the feasibility of the optimization problem. Section \ref{sec:simulation} shows some numerical experiments. Finally, in Section \ref{sec:conclusions} we draw the conclusions.

\subsection*{Notations}

In the following $\E$ denotes the mathematical expectation, $\Zbb$ the set of integers, $\Rbb$ the real line, and $\Cbb$ the complex plane. 
The symbol $\Hfrak_{n}$ represents the vector space (over the reals) of $n\times n$ Hermitian matrices, and $\Hfrak_{+,n}$ is the subset that contains positive definite matrices. 
The notation {$(\cdot)^{*}$ means taking complex conjugate transpose when applied to a matrix.
The symbol $\|\cdot\|$ may denote the norm of a matrix, a linear operator, or a function depending on the context.

\section{Problem formulation}
\label{sec:problem}

Suppose that we have a second-order stationary random field $\{\yb(\tb),\,\tb=(t_1,t_2,\dots,t_d)\in\Zbb^d\}$ where the positive integer $d$ is the dimension of the index set. For each $\tb\in\Zbb^d$, $\yb(\tb)$ is an $m$-dimensional zero mean complex random vector. The covariance is defined as $\Sigma_{\kb}:=\E\,\yb(\tb+\kb)\yb(\tb)^*$ which does not depend on $\tb$ by stationarity. In addition, we have the symmetry $\Sigma_{-\kb}=\Sigma_{\kb}^*$. The spectral density of the random field is defined as the Fourier transform of the matrix field $\{\Sigma_{\kb},\,\kb\in\Zbb^d\}$
\begin{equation}\label{Phi_spec_density}
\Phi(e^{i\thetab}):=\sum_{\kb\in\Zbb^d} \Sigma_{\kb}e^{-i\innerprod{\kb}{\thetab}},
\end{equation}
where $\thetab=(\theta_1,\theta_2,\dots,\theta_d)$ takes valued in $\Tbb^d:=(-\pi,\pi]^d$, $e^{i\thetab}$ is a shorthand for $(e^{i\theta_1},\dots,e^{i\theta_d})$, and $\innerprod{\kb}{\thetab}:=k_1\theta_1+\cdots+k_d\theta_d$ is the usual inner product in $\Rbb^d$. Given the symmetry of the covariances, one can easily verify that $\Phi$ is a Hermitian matrix-valued function on $\Tbb^d$.

Often in practice, a realization of the field $\yb$ is observed at a finite number of indices $\tb$ and we want to estimate the spectrum of the field from these observations. We shall proceed along the idea of \emph{rational covariance extension} \citep{Kalman}, starting by considering the covariances $\{\Sigma_{\kb},\,\kb\in\Lambda\}$, where $\Lambda\subset\Zbb^d$ is a finite index set such that
\begin{enumerate}
	\item $\zerob\in\Lambda$,
	\item $\kb\in\Lambda \implies -\kb\in\Lambda$.
\end{enumerate}
Then we aim to find a spectral density that matches these covariances. Formally, the problem is to find a function $\Phi:\Tbb^d \to \Hfrak_{+,m}$ that solves the integral equations
\begin{equation}\label{moment_eqns}
\int_{\Tbb^d}e^{i\innerprod{\kb}{\thetab}}\Phi(e^{i\thetab})\d m(\thetab) = \Sigma_{\kb},\text{ for all } \kb\in\Lambda
\end{equation}
given those $\Sigma_{\kb}$. Here
\begin{equation}
\d m(\thetab)=\frac{1}{(2\pi)^d}\prod_{j=1}^{d} \d\theta_j
\end{equation}
is the normalized Lebesgue measure over $\Tbb^d$.
The most common situation which will be the one referred to in our estimation procedure is the case in which $\Lambda$ is a cuboid centered at the origin.

When the integral equations \eqref{moment_eqns} are solvable, they usually have infinitely many solutions, and thus the problem above is not well-posed. The common approach in a still active line of research is to utilize entropy-like functionals as optimization criteria to select solutions.
In the same spirit of \citet{FMP-12}, we introduce the multidimensional version of the \emph{Itakura-Saito} (IS) distance between two bounded and coercive\footnote{A matrix spectral density $\Phi$ is bounded and coercive if there exist real numbers $M>\mu>0$ such that $\mu I_m \leq \Phi(e^{i\thetab}) \leq MI_m$ for all $\theta\in\Tbb^d$.} spectral densities
\begin{equation}\label{IS_dist_cont}
D(\Phi,\Psi):=\int_{\Tbb^d}\left(\log\det(\Phi^{-1}\Psi)+\trace[\Psi^{-1}(\Phi-\Psi)]\right)\d m.
\end{equation}
It is not difficult to see \citep[cf. e.g.,][p.~435]{LP15} that the latter is a pseudo-distance because $D(\Phi,\Psi) \geq 0$ and the equality holds if and only if $\Phi=\Psi$ (almost everywhere).

Our problem is now formulated as
\begin{equation}\label{primal_cont}
\underset{\Phi\in\Sfrak_m}{\text{minimize}}\ D(\Phi,\Psi) \quad \text{subject to } (\ref{moment_eqns}).
\end{equation}
Here the symbol $\Sfrak_m$ denotes the family of $\Hfrak_{+,m}$-valued functions defined on $\Tbb^d$ that are bounded and coercive. 
The spectral density function $\Psi$ which we call \emph{prior}, is given and it is interpreted as an extra piece of information that we have on the solution $\Phi$. More precisely, we want to find a solution to the moment equations \eqref{moment_eqns} that is closest to $\Psi$ as measured by $D$.
If no prior is available, we can select $\Psi=I$ in the spirit of maximizing the entropy, in which one aims to find the most unpredictable process with the prescribed moments. The maximum entropy (ME) paradigm is well accepted in the literature, as its meaningfulness has been discussed in \citet{csiszar1991least}. The unidimensional version of this optimization problem has been well studied in \citet{FMP-12}, which can be seen as a multivariate generalization of the scalar problem investigated in \citet{enqvist2008minimal}.

A similar covariance extension problem for random scalar fields has been studied in \citet{ringh2015multidimensional,RKL-16multidimensional,ringh2018multidimensional} using a different cost function (cf.~also \citet{KLR-16multidimensional} for a more general setting). 
Unlike the corresponding unidimensional problem, in the multidimensional case, the solution to the optimization problem, namely a \emph{spectral measure} that solves the moment equations, in general can contain a singular part.
This is a consequence of the fact that a certain integrability condition can fail to hold when the dimension $d\geq3$ (see \citet[Section 4]{RKL-16multidimensional} for details). However, such a singular measure is not unique, and its practical importance is so far still unclear.

An interesting exception is reported in \citet{ringh2015multidimensional} where a ``circulant'' version of the multidimensional covariance extension problem has been considered. There the spectral density function has support on a grid of $\Tbb^d$ (denoted as $\Tbb^d_\Nb$) and Fourier integrals such as those in \eqref{moment_eqns} are replaced by (inverse) discrete Fourier transforms. It is shown in that paper that given a \emph{positive} trigonometric polynomial $P$ on $\Tbb^d_\Nb$, there exists a unique polynomial $Q$ that is also positive on $\Tbb^d_\Nb$, such that the rational function $\Phi=P/Q$ solves the moment equations. In this case, no singular measure arises, which is exactly analogous to the main result in \citet{LPcirculant-13} that treats the unidimensional problem.

Next we shall mainly work on the discrete version of the optimization problem in \eqref{primal_cont}. Let us set up the notation first. The product set $\Tbb^d=(-\pi,\pi]^d$ is discretized such that $N_j$ equidistant points are selected from the $j$-th factor $(-\pi,\pi]$. More precisely, let us fix
\begin{equation}\label{Nb}
\Nb=(N_1,N_2,\dots,N_d),
\end{equation}
and define the finite index set (as a subset of $\Zbb^d$)
\begin{equation}\label{Zd_N}
\Zbb^d_\Nb:=\left\{ \lb=(\ell_1,\dots,\ell_d) : 0\leq \ell_j \leq N_j-1,\, j=1,\dots,d \right\}.
\end{equation}
The set $\Zbb_\Nb^d$ has cardinality $|\Nb|:=\prod_{j=1}^{d}N_j$.
The discretization of $\Tbb^d$ can then be expressed as
\begin{equation*}
\Tbb^d_\Nb:=\left\{ \left( \frac{2\pi}{N_1}\ell_1,\dots,\frac{2\pi}{N_d}\ell_d \right) : \lb\in\Zbb_\Nb^d \right\}.
\end{equation*}
Moreover, let $\zetab_{\lb}:=\left(\zeta_{\ell_1},\dots,\zeta_{\ell_d}\right)$ be an element of the discretized $d$-torus with $\zeta_{\ell_j}=e^{i2\pi\ell_j/N_j}$ and define $\zetab^{\kb}_{\lb}:=\prod_{j=1}^{d}\zeta_{\ell_j}^{k_j}$.
Define next a discrete measure with equal mass on the grid points in $\Tbb^d_\Nb$:
\begin{equation}
\d\nu_\Nb(\thetab) = \sum_{\lb\in\Zbb^d_\Nb} \delta ( \theta_1-\frac{2\pi}{N_1}\ell_1,\dots,\theta_d-\frac{2\pi}{N_d}\ell_d ) \prod_{j=1}^{d}\frac{\d\theta_j}{N_j}.
\end{equation}
The IS distance between spectra defined on $\Tbb^d_\Nb$ takes the form 
{\small \begin{align} D_\Nb(\Phi,\Psi):=\int_{\Tbb^d}\left(\log\det(\Phi^{-1}\Psi)+\trace[\Psi^{-1}(\Phi-\Psi)]\right)\d\nu_\Nb .\nonumber
	\end{align}}
Our discretized optimization problem can be written as
\begin{subequations}\label{primal_discrete}
	\begin{align}
	& \underset{\Phi(\zetab_{\lb})>0,\,\forall\,\lb\in\Zbb^d_\Nb}{\text{minimize}}
	D_\Nb(\Phi,\Psi) \label{cost_discrete} \\
	& \text{s.t. }
	\int_{\Tbb^d}e^{i\innerprod{\kb}{\thetab}}\Phi(e^{i\thetab})\d\nu_\Nb(\thetab) = \Sigma_{\kb},\text{ for all } \kb\in\Lambda \label{constraints_discrete}
	\end{align}
\end{subequations}
where the integrals
\begin{equation}
\int_{\Tbb^d} f(e^{i\thetab})\d\nu_\Nb(\thetab)=\frac{1}{|\Nb|} \sum_{\lb\in\Zbb^d_\Nb} f(\zetab_{\lb})
\end{equation}
are essentially Riemann sums. 

There are two reasons to prefer the discrete formulation \eqref{primal_discrete}.
\begin{enumerate}[i)]
	\item From the numerical aspect, we will have to discretize the problem when implementing an algorithm on a computer, and we may as well treat the discretized problem in the first place. Moreover, the fast Fourier transform (FFT) can be used to compute quantities such as the moments.
	\item The number of available data is anyway finite so that a discrete theory appears natural. As we will see, it will provide nice and elegant theory with no need to consider {\em ad hoc} singular measures.
\end{enumerate}
In addition, the discrete spectrum defined over $\Tbb_\Nb^d$ has a probabilistic interpretation of corresponding to a \emph{periodic} stationary field, as we shall explain below.

\subsection{Spectral representation of periodic stationary fields}\label{subsec:spec_rep}

Let $\yb(\tb)$ be a second-order stationary random complex $m$-vector field defined over $\Zbb^d$. By stationarity, the field admits a representation \citep[cf.][]{yaglom1957some}\footnote{The author of \citet{yaglom1957some} used the attributive ``homogeneous'' in place of ``stationary''.}
\begin{equation}\label{y_t_spec}
\yb(\tb) = \int_{\Tbb^d} e^{i\innerprod{\tb}{\thetab}} \d \hat{\yb}(\thetab),\quad \tb\in\Zbb^d,
\end{equation}
where $\hat{\yb}(\thetab)=[\,\hat{y}_1(\thetab),\dots, \hat{y}_m(\thetab)\,]$ is a vector of random measures on $\Tbb^d$.

Next, let us impose the following periodicity assumption. Suppose that for any $\tb=(t_1,\dots,t_d)$,
\begin{equation}
\begin{split}
\yb(\tb) & = \yb(t_1+N_1,t_2,t_3,\dots,t_d) \\
& \vdots \\
\yb(\tb) & = \yb(t_1,\dots,t_{j-1},t_j+N_j,t_{j+1},\dots,t_d) \\
& \vdots \\
\yb(\tb) & = \yb(t_1,t_2,\dots,t_{d-1},t_d+N_d)
\end{split}
\end{equation}
almost surely. In other words, the field $\yb(\tb)$ is periodic with a period of $N_j$ (a positive integer) in the $j$-th dimension. Using the spectral representation \eqref{y_t_spec}, the periodicity assumption in $t_1$ (the first dimension) implies
\begin{equation}
\int_{\Tbb^d} e^{i\innerprod{\tb}{\thetab}} (e^{i N_1 \theta_1} - 1) \, \d\hat{\yb}(\thetab) = \zerob \quad \text{a.s.}
\end{equation}
Multiply both sides of the equation with their complex conjugate transposes, take the expectation, and we get
\begin{equation}\label{spec_mat}
\int_{\Tbb^d} |e^{i N_1 \theta_1} - 1|^2 \, \d\Fb(\thetab) = \zerob
\end{equation}
where $\d\Fb(\thetab) = \E\, \d\hat{\yb}(\thetab) \d\hat{\yb}(\thetab)^*$ and $\Fb(\cdot)$, called the \emph{spectral matrix} of the field $\yb(\tb)$, is a Hermitian nonnegative definite matrix of complex measures defined on the Borel subsets of $\Tbb^d$. The equality \eqref{spec_mat} implies that the support of $\d\Fb$ must be contained in $\Tbb_{N_1} \times \Tbb^{d-1}$ where
\begin{equation}
\Tbb_{N_1}:= \left\{ \frac{2\pi}{N_1}\ell_1 \,:\, 0\leq\ell_1\leq N_1-1 \right\}.
\end{equation}
Repeat the argument in each dimension, and we conclude that the support of $\d\Fb$ is in fact contained in the grid $\Tbb_\Nb^d$.

Apparently, the random field being periodic leads to the periodicity of the covariance field $\{\Sigma_{\kb}\}_{\kb\in\Zbb^d}$. More precisely, we have that in the $j$-th dimension
\begin{equation}
\Sigma_{\kb} = \Sigma_{(k_1,\dots,k_{j-1},k_j+N_j,k_{j+1},\dots,k_d)}.
\end{equation}
Therefore, we can restrict our attention to one particular $d$-cuboid $\{\Sigma_{\kb}\}_{\kb\in\Zbb_\Nb^d}$. Moreover, combining the periodicity with the usual symmetry $\Sigma_{-\kb}=\Sigma_{\kb}^*$ induced by the stationarity assumption, we have 
\begin{equation}
\Sigma_{(-k_1,\dots,-k_{j-1},-k_j+N_j,-k_{j+1},\dots,-k_d)} = \Sigma_{\kb}^*,
\end{equation}
which is another kind of symmetry for the covariances around the ``end points'' of the $d$-cuboid.
In the unidimensional case ($d=1$), the above equality reduces to $\Sigma_{-k+N} = \Sigma_k^*$, which adds a block-circulant structure on the covariance matrix $\Sigmab:= \E\,\Yb \Yb^*$, where $\Yb$ is a long column vector of random variables obtained by stacking $\yb(0),\dots,\yb(N)$ together \citep[see e.g.,][]{CFPP-11}. In the $2$-d case for example, suppose that we know the covariances at the indices $\{(0,0),(1,0),(0,1),(1,1)\}$. Then automatically, we have
\begin{equation}
\begin{split}
\Sigma_{(N_1-1,0)} & = \Sigma_{(-1,0)} = \Sigma_{(1,0)}^* \, , \\
\Sigma_{(0,N_2-1)} & = \Sigma_{(0,-1)} = \Sigma_{(0,1)}^* \, , \\
\Sigma_{(N_1-1,N_2-1)} & = \Sigma_{(-1,-1)} = \Sigma_{(1,1)}^* \, .
\end{split}
\end{equation}

\section{The dual optimization problem}\label{sec:dual}

In this section, we will elaborate how to solve the optimization problem \eqref{primal_discrete} via duality. 
In the first place, we make a feasibility assumption.
\begin{assumption}[Feasibility]\label{assump_feasibility}
	There exists a function $\Phi: \Tbb^d_\Nb \to \Hfrak_{+,m}$ such that
	the constraints \eqref{constraints_discrete} hold with the given $\{\Sigma_\kb\}_{\kb\in\Lambda}$.
\end{assumption}

\begin{remark}
	In \citet{ringh2015multidimensional}, the feasibility assumption is stated in terms of a dual cone formulation. See also \citet{KLR-16multidimensional,RKL-16multidimensional,ringh2018multidimensional}. As we shall see next, the feasibility assumption plays an important role in the development of the theory, and later in Section~\ref{sec:spec_est}, an estimator of the covariances will be proposed such that the above assumption is satisfied in practice. Another way to deal with the feasibility can be found in \citet{enqvist2007approximative,ringh2018multidimensional}.
\end{remark}

Notice that due to discretization, the optimization problem \eqref{primal_discrete} is a finite dimensional problem subject to a finite number of linear equality constraints. In the literature of optimization theory, there are methods to handle this type of problems directly. However, it is usually more convenient to work with the dual problem, because the number of dual variables is proportional to the cardinality of the set $\Lambda$ while the number of primal variables is proportional to the grid size $|\Nb|$ which is usually much larger than the former.

The discretized IS-distance can be rewritten as
\begin{align*}
&D_\Nb(\Phi,\Psi)=\\
&\qquad \int_{\Tbb^d}\hspace{-0.2cm} \left( -\log\det\Phi +\trace(\Psi^{-1}\Phi) +\log\det\Psi  \right)\d\nu_\Nb-m.
\end{align*}

One can see that the last two terms in the integral do not depend on $\Phi$, and thus can be neglected from the cost function. Let us form the Lagrangian
\begin{equation*}
\begin{split}
\Lcal_\Psi( &\Phi,\Qb)   = \int_{\Tbb^d} \left( -\log\det\Phi +\trace(\Psi^{-1}\Phi) \right)\d\nu_\Nb \\
&  +\sum_{\kb\in\Lambda} \trace \left[ Q_\kb \left(\int_{\Tbb^d} e^{i\innerprod{\kb}{\thetab}} \Phi(e^{i\thetab}) \d\nu_\Nb(\thetab) - \Sigma_{\kb}\right)^* \right] \\
& = \int_{\Tbb^d} \left\{ -\log\det\Phi +\trace[(\Psi^{-1}+Q)\Phi] \right\}\d\nu_\Nb-\innerprod{\Qb}{\Sigmab}
\end{split}
\end{equation*}
where, the variable $\Qb=\{Q_\kb\}_{\kb\in\Lambda}$ contains the Lagrange multipliers such that each $Q_\kb\in\Cbb^{m \times m}$ and $Q_\zerob$ is Hermitian, $Q(e^{i\thetab}):=\sum_{\kb\in\Lambda} Q_\kb e^{-i\innerprod{\kb}{\thetab}}$ is a matrix trigonometric polynomial of several variables, $\Sigmab=\{\Sigma_\kb\}_{\kb\in\Lambda}$ consists of the covariance data, and $\innerprod{\Qb}{\Sigmab}:=\sum_{\kb\in\Lambda}\trace(Q_\kb \Sigma_{\kb}^*)$.

For a fixed $\Qb$, consider the problem
\begin{equation*}
\underset{\Phi(\zetab_{\lb})>0,\,\forall\,\lb\in\Zbb^d_\Nb}{\text{minimize}}\ \Lcal_\Psi(\Phi,\Qb).
\end{equation*}
The function $\Lcal_\Psi(\cdot,\Qb)$ is strictly convex in the feasible set. The directional derivative of the Lagrangian in the direction $\delta\Phi: \Tbb^d_\Nb \to \Hfrak_{m}$ can be computed as
\begin{equation*}
\begin{split}
\delta \Lcal_\Psi(\Phi,\Qb;\delta\Phi) & = \int_{\Tbb^d} \trace[(-\Phi^{-1}+\Psi^{-1}+Q) \delta\Phi] \d\nu_\Nb  \\
& = \int_{\Tbb^d}\innerprod{-\Phi^{-1}+\Psi^{-1}+Q}{\delta\Phi} \d\nu_\Nb \\
\end{split}
\end{equation*}
where we have used the fact that the directional derivative of $\log\det X$ for $X>0$ is given by
\begin{equation*}
\delta \log \det (X;\delta X) = \trace (X^{-1} \delta X).
\end{equation*}
We impose the first variation $\delta \Lcal_\Psi(\Phi,\Qb;\delta\Phi)$ to vanish in any direction $\delta\Phi$. In particular, taking $\delta\Phi=-\Phi^{-1}+\Psi^{-1}+Q$ implies
\begin{equation*}
\Phi=\Phi^{\circ}(\Qb):=(\Psi^{-1}+Q)^{-1},\quad \forall \thetab \in \Tbb^d_\Nb.
\end{equation*}

Since $\Phi$ is required to be positive definite on the grid $\Tbb^d_\Nb$, the Lagrange multiplier $\Qb$ must be constrained to the set
\begin{equation*}
\Lscr_+:=\left\{ \{Q_\kb\}_{\kb \in\Lambda} : (\Psi^{-1}+Q)(\zetab_{\lb})>0,\ \forall \lb\in\Zbb^d_\Nb \right\}.
\end{equation*}
By the continuous dependence of eigenvalues on the matrix entries, one can verify that $\Lscr_+$ is an open set. Insert $\Phi^\circ$ into the Lagrangian to yield the dual problem to maximize the expression
\begin{equation*}
\int_{\Tbb^d} \log\det(\Psi^{-1}+Q) \d\nu_\Nb-\innerprod{\Qb}{\Sigmab}+m.
\end{equation*}
Hereafter we will instead consider the equivalent problem
\begin{equation}\label{J_dual}
\underset{\Qb\in\Lscr_+}{\text{minimize}} \quad J_\Psi(\Qb):=\langle\Qb,\Sigmab\rangle-\int_{\Tbb^d}\log\det(\Psi^{-1}+Q)\d\nu_\Nb,
\end{equation}
and call $J_\Psi$ the dual function. The main theorem of this section is stated below, and the proof will be presented subsequently.

\begin{theorem}\label{thm_main}
	If Assumption $\ref{assump_feasibility}$ holds, then the dual function $J_\Psi$ has a unique minimizer $\Qb^\circ\in\Lscr_+\,$. Moreover, the spectral density $\Phi^{\circ}(\Qb^\circ)=(\Psi^{-1}+Q^\circ)^{-1}$ defined over $\Tbb_\Nb^d$ solves the discretized primal problem $\eqref{primal_discrete}$.
\end{theorem}

\subsection{Uniqueness of the minimizer}\label{subsec:uniqueness}

We claim that if a minimizer $\Qb^\circ$ of $J_\Psi$ \emph{exists}, it is unique. This is a consequence of the strict convexity of the cost function. To see this, let us compute the first and second variations of $J_\Psi$.
\begin{align*}
& \delta J_\Psi (\Qb;\delta\Qb) \nonumber\\
& = \innerprod{\delta\Qb}{\Sigmab} - \int_{\Tbb^d} \trace \left[ (\Psi^{-1}+Q)^{-1} \sum_{\kb\in\Lambda} \delta Q_\kb e^{-i\innerprod{\kb}{\thetab}} \right] \d\nu_\Nb \nonumber \\
& = \trace \left\{ \sum_{\kb\in\Lambda} \delta Q_\kb \left[ \Sigma^*_\kb - \int_{\Tbb^d} e^{-i\innerprod{\kb}{\thetab}}(\Psi^{-1}+Q)^{-1} \d\nu_\Nb \right] \right\}.\nonumber 
\end{align*}
It is not difficult to see that the differential $\delta J_\Psi(\Qb;\delta\Qb)$ is continuous in $\Qb$ for any fixed direction $\delta\Qb$. In fact, it amounts to showing the continuity of the term
\begin{equation*}
\begin{split}
& \int_{\Tbb^d} e^{-i\innerprod{\kb}{\thetab}}(\Psi^{-1}+Q)^{-1} \d\nu_\Nb \\
= & \frac{1}{|\Nb|} \sum_{\lb\in\Zbb^d_\Nb} \zetab_{\lb}^{\kb} \left( \Psi^{-1}(\zetab_{\lb}) + Q(\zetab_{\lb}) \right)^{-1},
\end{split}
\end{equation*}
which is trivial because of the finite summation. Furthermore, due to the smoothness of the matrix inversion map $X\mapsto X^{-1}$, the function $J_\Psi$ is smooth over $\Lscr_+$.

A feasible $\Qb^\circ\in\Lscr_+$ that annihilates directional derivatives in every direction $\delta\Qb$ must satisfy the relation
\begin{equation}\label{stationary_cond}
\int_{\Tbb^d} e^{i\innerprod{\kb}{\thetab}}(\Psi^{-1}+Q^{\circ})^{-1} \d\nu_\Nb = \Sigma_\kb,\quad\forall \kb\in\Lambda.
\end{equation}
In other words, the spectral density $\Phi^\circ(\Qb^\circ)=(\Psi^{-1}+Q^\circ)^{-1}$ is a solution to the discretized moment equations.

The second-order derivative (differential) of $J_\Psi$ at $\Qb$ is
\begin{equation*}
\begin{split}
&\delta^2 J_\Psi(\Qb;\delta\Qb^{(1)},\delta\Qb^{(2)}) \\ & = \trace\Bigg[ \sum_{\kb\in\Lambda}\sum_{\tilde \kb\in\Lambda} \delta Q^{(1)}_\kb \int_{\Tbb^d} e^{-i\innerprod{\kb+\tilde\kb}{\thetab}} (\Psi^{-1}+Q)^{-1}  \\ & \hspace{1.2cm} \times \delta Q^{(2)}_{\tilde\kb}  (\Psi^{-1}+Q)^{-1} \d\nu_\Nb \Bigg] \\
& = \trace \int_{\Tbb^d} \delta Q^{(1)} (\Psi^{-1}+Q)^{-1} \delta Q^{(2)}(\Psi^{-1}+Q)^{-1} \d\nu_\Nb
\end{split}
\end{equation*}
which is understood as a bilinear function in $\delta\Qb^{(k)}$, $k=1,2$. In the above computation, we have used the fact that the differential of the map $X \mapsto X^{-1}$ at $X$ is given by $\delta X \mapsto -X^{-1} \delta X X^{-1}$.

For $\Qb\in\Lscr_+$, write $\Phi^\circ=(\Psi^{-1}+Q)^{-1}$ for short. Since $\Phi^\circ(\zetab_{\lb})>0$ for all $\lb\in\Zbb_\Nb^d$, we can perform the Cholesky factorization $\Phi^\circ = LL^*$ where each quantity here depends on the discrete frequency $\zetab_{\lb}$.
Therefore we have the second variation
\begin{equation*}
\begin{split}
\delta^2& J_\Psi(\Qb;\delta\Qb,\delta\Qb)   = \trace \int_{\Tbb^d} \delta Q \, \Phi^\circ \, \delta Q \, \Phi^\circ \d\nu_\Nb \\
& = \frac{1}{|\Nb|} \trace \sum_{\lb\in\Zbb^d_\Nb} \delta Q(\zetab_{\lb}) \, \Phi^\circ(\zetab_{\lb}) \, \delta Q(\zetab_{\lb}) \, \Phi^\circ(\zetab_{\lb}) \\
& = \frac{1}{|\Nb|} \trace \sum_{\lb\in\Zbb^d_\Nb} L^*(\zetab_{\lb}) \, \delta Q(\zetab_{\lb}) \, \Phi^\circ(\zetab_{\lb}) \, \delta Q(\zetab_{\lb}) \, L(\zetab_{\lb}) \geq 0.
\end{split}
\end{equation*}
Since $\Phi^\circ(\zetab_{\lb})$ is positive definite for any $\zetab_{\lb}$, the second variation of $J_\Psi$ is equal to zero if and only if the polynomial $\delta Q(\zb)$ vanishes identically on the discrete $d$-torus corresponding to the frequencies in $\Tbb^d_\Nb$. We shall make the following innocuous assumption for the strict positivity of the second variation. For $j=1,\dots,d$, let $n_j:=\max\{|k_j| : \kb\in\Lambda\}$, the largest index in the $j$-th dimension.

\begin{assumption}\label{assump_N}
	The positive integers used to define the grid $\Tbb_\Nb^d$ in \eqref{Nb} satisfy $N_j>2n_j$ for $j=1,\dots,d$.
\end{assumption}
If Assumption \ref{assump_N} holds, then according to Lemma~1 in \citet{ringh2015multidimensional}, $\delta^2 J_\Psi(\Qb;\delta\Qb,\delta\Qb)=0$ implies $\delta Q(\zb) \equiv 0$, i.e. the zero polynomial. We can now conclude that the dual problem \eqref{J_dual} is strictly convex. The uniqueness of the solution is an easy consequence of the strict convexity, provided that a solution exists. The existence question is indeed the difficult part, which is the content of the next subsection.

\subsection{Existence of an interior minimizer}

In this subsection, we will show that a minimizer $\Qb^\circ$ of $J_\Psi$ exists in the open set $\Lscr_+$. Let us begin by defining the boundary of the set $\Lscr_+$,
\begin{align*}
\partial\Lscr_+ := \{ \Qb & : \Psi^{-1}+Q\geq 0 \text{ on } \Tbb^d_\Nb\nonumber \\  &   \text{    and is singular at least in one point} \} \end{align*}
and the closure $\overline{\Lscr_+}:=\Lscr_+\cup\partial\Lscr_+$.

\begin{proposition}\label{prop_J_bounded_below}
	For any $\Qb\in\Lscr_+$, the function value $J_\Psi(\Qb)$ is bounded from below.
\end{proposition}

In particular, this implies that the dual function $J_\Psi$ cannot attain a value of $-\infty$ on $\Lscr_+$.
To prove the proposition, we need the following lemma.

\begin{lemma}\label{lem_innerprod_inequal}
	If Assumption $\ref{assump_feasibility}$ holds, then there exist two real constants $\mu>0$ and $\alpha$ such that for any $\Qb\in\overline{\Lscr_+}\,$,
	\begin{equation}\label{innerprod_inequal}
	\innerprod{\Qb}{\Sigmab} \geq \alpha + \mu \trace \int_{\Tbb^d} (\Psi^{-1}+Q) \d\nu_\Nb.
	\end{equation}
\end{lemma}

\begin{pf}
	The feasibility assumption says that there exists a function $\Phi_0$ that is positive definite on $\Tbb_\Nb^d$ and satisfies the moment equations in \eqref{constraints_discrete} for all $\kb\in\Lambda$. One can see that there exists a positive constant $\mu$ such that $\Phi_0(\zetab_{\lb}) \geq \mu I$ for all $\zetab_{\lb}$ simply because $\Tbb_\Nb^d$ contains a finite number of elements. We have
	\begin{equation}\label{innerprod_inequal_deriv}
	\begin{split}
	\innerprod{\Qb}{\Sigmab}& := \sum_{\kb\in\Lambda} \trace (Q_\kb\Sigma_\kb^*) \\ & = \sum_{\kb\in\Lambda} \trace \left( Q_\kb \int_{\Tbb^d} e^{-i\innerprod{\kb}{\thetab}} \Phi_0(e^{i\thetab}) \d\nu_\Nb \right) \\
	& = \trace \int_{\Tbb^d} Q \, \Phi_0 \, \d\nu_\Nb \\
	& = \trace \int_{\Tbb^d} (\Psi^{-1}+Q) \Phi_0 \, \d\nu_\Nb \underbrace{- \trace \int_{\Tbb^d} \Psi^{-1} \, \Phi_0 \, \d\nu_\Nb}_{:=\alpha} \\
	& = \alpha + \trace \int_{\Tbb^d} L_1^* \, \Phi_0 \, L_1 \, \d\nu_\Nb \\
	& \geq \alpha + \mu \trace \int_{\Tbb^d} (\Psi^{-1}+Q)\d\nu_\Nb,
	\end{split}
	\end{equation}
	where $L_1$ (as a function of $\zetab_{\lb}$) is the Cholesky factor of $\Psi^{-1}+Q\geq 0$ since $\Qb\in\overline{\Lscr_+}$, and the constant $\alpha<0$.
\end{pf}

\begin{pf}[Proof of Proposition $\ref{prop_J_bounded_below}$]
	By Lemma \ref{lem_innerprod_inequal}, we have
	\begin{equation}
	\begin{split}
	J_\Psi(\Qb) & = \innerprod{\Qb}{\Sigmab}-\int_{\Tbb^d} \log \det (\Psi^{-1}+Q) \d\nu_\Nb \\
	& \geq \alpha + \trace \int_{\Tbb^d} \left[ \mu(\Psi^{-1}+Q) - \log (\Psi^{-1}+Q) \right] \d\nu_\Nb,
	\end{split}
	\end{equation}
	where we have used the fact that $\log \det X = \trace \log X$ for $X>0$. The matrix $(\Psi^{-1}+Q)(\zetab_{\lb})$ is positive definite since we have the condition $\Qb\in\Lscr_+$, and let $\{\lambda_j(\zetab_{\lb})\}_{j=1}^m$ be its eigenvalues. Then the above inequality can be written as
	\begin{equation}
	\Jbb_\Psi(\Qb) \geq \alpha + \int_{\Tbb^d} \sum_{j=1}^{m} ( \mu\lambda_j - \log \lambda_j) \, \d\nu_\Nb.
	\end{equation}
	Define the function $\rho(\lambda_1,\dots,\lambda_m):=\sum_{j=1}^{m} ( \mu\lambda_j - \log \lambda_j)$. One can easily verify that $\rho$ is strictly convex in the orthant $\{ (\lambda_1,\dots,\lambda_m) : \lambda_j>0,\ \forall j=1,\dots,m \}$. Its minimizer is the stationary point $\lambda_j=1/\mu>0$, $j=1,\dots,m$ with a minimum value $m(1+\log \mu)$. We thus obtain a uniform lower bound for the function value $J_\Psi(\Qb)$, namely
	\begin{equation}
	J_\Psi(\Qb) \geq \alpha + m(1+\log \mu).
	\end{equation}
\end{pf}

Following the same reasoning as above, one can show without difficulty that if a sequence $\{\Qb_k\}_{k\geq 1}$ in $\Lscr_+$ tends to some $\bar{\Qb}\in \partial\Lscr_+$, then
\begin{equation}\label{J_at_boundary}
\lim_{k\to\infty}J_\Psi(\Qb_k) = +\infty.
\end{equation}
This is true because the usual Lebesgue measure is replaced by the discrete measure $\d\nu_\Nb$ in the integrals during the discretization of the optimization problem. In sharp contrast, as reported in \citet{FMP-12} for the unidimensional problem, the dual function may still have a finite value on the boundary of the feasible set when the integration is done with respect to the Lebesgue measure.

Take a sufficiently large real number $r$. Any point $\Qb$ satisfying $J_\Psi(\Qb)\leq r$ must be away from the boundary $\partial\Lscr_+$ due to \eqref{J_at_boundary}. Therefore, we can define the (nonempty) sublevel set of $J_\Psi$ as a subset of $\Lscr_+$:
\begin{equation}
J_\Psi^{-1}(-\infty,r] := \{ \Qb\in\Lscr_+ : J_\Psi(\Qb)\leq r \}.
\end{equation}
Moreover, the sublevel set is closed because of the continuity of $J_\Psi$. 

Define
\begin{equation}\label{Q_norm}
\|\Qb\|:= \sqrt{ \sum_{\kb\in\Lambda} \trace (Q_\kb Q^*_\kb) },
\end{equation}
the norm of the object $\Qb=\{Q_\kb : \kb\in\Lambda \}$ that contains the Lagrange multipliers. We have the next lemma.
\begin{lemma}\label{lem_inf_Q}
	Let $\{\Qb_k\}_{k\geq 1}$ be a sequence in ${\Lscr_+}$ such that $\|\Qb_k\|\to\infty$ as $k\to \infty$. Then there exists a subsequence $\{\Qb_{k_j}\}_{j\geq 1}$ such that
	\begin{equation}\label{limit_inf_Q}
	\lim_{j\to\infty} J_\Psi(\Qb_{k_j}) = +\infty.
	\end{equation}
\end{lemma}

\begin{pf}
	Given the sequence $\{\Qb_k\}$, define $\Qb_k^0:=\Qb_k / \|\Qb_k\|$, which necessarily implies that $Q_k^0(e^{i\thetab})=Q_k(e^{i\thetab}) / \|\Qb_k\|$. Moreover, since each $\Qb_k\in {\Lscr_+}$, we have $\Psi^{-1}+Q_k > 0$ on $\Tbb_\Nb^d$. Consequently, the function
	\begin{equation}
	\Psi^{-1} + Q_k^0 = \frac{1}{\|\Qb_k\|} (\Psi^{-1}+Q_k) + \left( 1 - \frac{1}{\|\Qb_k\|} \right) \Psi^{-1}
	\end{equation}
	is positive definite on $\Tbb_\Nb^d$ for all sufficiently large $k$ since $\|\Qb_k\|\to\infty$. To summarize, the sequence $\{\Qb_k^0\}$ lives on the unit surface $\|\Qb\|=1$ (a compact set due to finite dimensionality), and we have $\Qb_k^0\in\Lscr_+$ for $k$ large enough.

	From Lemma \ref{lem_innerprod_inequal}, we know that $\innerprod{\Qb}{\Sigmab}\geq \alpha$ for any $\Qb\in\overline{\Lscr_+}$, because the second term on the right side of \eqref{innerprod_inequal} is nonnegative.
	Hence
	\begin{equation}
	\innerprod{\Qb_k^0}{\Sigmab} = \frac{1}{\|\Qb_k\|} \innerprod{\Qb_k}{\Sigmab} \geq \frac{\alpha}{\|\Qb_k\|} \to 0.
	\end{equation}
	Define the real quantity $\eta:= \liminf_{k\to\infty} \innerprod{\Qb_k^0}{\Sigmab}$. Then it must hold that $\eta\geq 0$. By a property of the limit inferior, we know that $\{\Qb_k^0\}$ has a subsequence $\{\Qb^0_{k_\ell}\}$ such that $\innerprod{\Qb_{k_\ell}^0}{\Sigmab}\to\eta$ as $\ell\to\infty$. Since $\{\Qb^0_{k_\ell}\}_{\ell\geq 1}$ is contained on the unit surface, it has a convergent subsequence denoted by $\{\Qb^0_{k_j}\}_{j\geq 1}$. Define the limit
	\begin{equation}
	\Qb^0_{\infty}:=\lim_{j\to\infty}\Qb_{k_j}^0.
	\end{equation}
	Then by the continuity of the inner product, we have $\eta=\innerprod{\Qb^0_{\infty}}{\Sigmab}$.
	
	Next, we show that $\Qb^0_\infty \in \Lscr_+$. Since $\Qb_k\in {\Lscr_+}$, it holds that $\Psi^{-1}+Q_k>0$ on $\Tbb_\Nb^d$ for all $k$. This implies that
	\begin{equation}
	\frac{\Psi^{-1}}{\|\Qb_{k_j}\|}+Q^0_{k_j} > 0 \text{ on } \Tbb_\Nb^d, \ \forall j.
	\end{equation}
	The function on the left side of the above inequality converges uniformly to the polynomial $Q^0_\infty$ on $\Tbb_\Nb^d$. Hence we must have $Q^0_\infty\geq 0$ on $\Tbb_\Nb^d$. As a consequence, $\Psi^{-1}+Q^0_\infty>0$ on $\Tbb_\Nb^d$ and indeed $\Qb^0_\infty \in \Lscr_+$.
	
	The next step is to prove that $\eta=\innerprod{\Qb^0_{\infty}}{\Sigmab}>0$. Following the computation in \eqref{innerprod_inequal_deriv}, we arrive at
	\begin{equation}
	\begin{split}
	\innerprod{\Qb^0_\infty}{\Sigmab} & = \trace \int_{\Tbb^d} Q^0_\infty \, \Phi_0 \, \d\nu_\Nb \\
	& = \trace \int_{\Tbb^d} L_2^* \, Q^0_\infty \, L_2 \, \d\nu_\Nb,
	\end{split}
	\end{equation}
	where $\Phi_0 = L_2 L_2^*$ is the point-wise Cholesky factorization. Since we have just proved that $Q^0_\infty$ is positive semidefinite on $\Tbb_\Nb^d$, the same is true for the function $L_2^* Q^0_\infty L_2$. Thus, $\innerprod{\Qb^0_{\infty}}{\Sigmab}=0$ implies that the polynomial $Q^0_\infty(\zb)$ vanishes identically on the discrete $d$-torus. By \citet[Lemma~1]{ringh2015multidimensional}, we must have $\Qb^0_\infty=\zerob$, which is a contradiction since we also have $\|\Qb^0_\infty\|=1$. Therefore, it must hold that $\eta>0$.
	
	Finally, since $\innerprod{\Qb^0_{k_j}}{\Sigmab} \to \eta>0$ as $j\to\infty$, there exists an integer $K>0$ such that $\forall j>K$, $\innerprod{\Qb^0_{k_j}}{\Sigmab} > \eta/2$.
	{\small \begin{equation}
		\begin{split}
		& \lim_{j\to\infty} J_\Psi(\Qb_{k_j})   = \lim_{j\to\infty} \innerprod{\Qb_{k_j}}{\Sigmab} - \int_{\Tbb^d} \log \det (\Psi^{-1} + Q_{k_j}) \d\nu_\Nb \\
		& = \lim_{j\to\infty} \|\Qb_{k_j}\| \innerprod{\Qb_{k_j}^0}{\Sigmab}\\ & \hspace{0.8cm}- \int_{\Tbb^d} \log \det \|\Qb_{k_j}\| \left( \frac{\Psi^{-1}}{\|\Qb_{k_j}\|} + Q^0_{k_j} \right) \d\nu_\Nb \\
		& \geq \lim_{j\to\infty} \frac{\eta}{2} \|\Qb_{k_j}\| - m \log \|\Qb_{k_j}\|   \\ & \hspace{0.8cm}-\int_{\Tbb^d} \log \det \left( \frac{\Psi^{-1}}{\|\Qb_{k_j}\|} + Q^0_{k_j} \right) \d\nu_\Nb
		\end{split}
		\end{equation}}
	The function $\Psi^{-1} / \|\Qb_{k_j}\| + Q^0_{k_j}$ has bounded norm over $\Tbb_\Nb^d$. Hence the integral of its $\log \det$ is bounded from above. Comparing linear and logarithmic growth, we can make the conclusion \eqref{limit_inf_Q}.
\end{pf}

As a consequence of Lemma \ref{lem_inf_Q}, the sublevel set $J_\Psi^{-1}(-\infty,r]$ has to be bounded. Therefore, it is a compact subset of $\Lscr_+$. By the extreme value theorem, the function $J_\Psi$ has a minimum over the sublevel set, and the minimizer $\Qb^\circ$ is in $\Lscr_+$. This concludes the existence proof.

Now we have shown that the dual problem \eqref{J_dual} has a unique solution. The remaining claim of Theorem \ref{thm_main} follows from strong duality, i.e., zero duality gap.

\section{Well-posedness}\label{sec:well-posed}

In the previous section, we have demonstrated that the optimization problem \eqref{primal_discrete} has a unique solution via duality. 
These properties, although necessary, are far from being sufficient for the problem to make sense from the engineering point of view. Indeed, to this end it is fundamental that the solution depends continuously on the problem data. In this section we show that this is the case and that, in fact,   
the problem is well-posed in the sense of Hadamard. To show this, we are left to establish the smooth dependence of the solution $\Qb$ on the problem data, namely the prior spectral density $\Psi$ and the covariance matrices $\Sigmab=\{\Sigma_{\kb}\}_{\kb\in\Lambda}$. The argument is built upon the classical inverse and implicit function theorems.

Define first the linear operator that sends a Hermitian-matrix-valued function on $\Tbb_\Nb^d$ to its discrete Fourier coefficients with indices in the set $\Lambda$
\begin{equation}
\Gamma:\, \Phi \mapsto \left\{ \Sigma_{\kb}=\int_{\Tbb^d} e^{i\innerprod{\kb}{\thetab}} \, \Phi \, \d\nu_\Nb \right\}_{\kb\in\Lambda}.
\end{equation}
Let
\begin{equation}
\Mscr_+:=\left\{ \Sigmab = \Gamma(\Phi) \,:\, \Phi(\zetab_{\lb})>0, \ \forall \lb\in\Zbb_\Nb^d \right\}
\end{equation}
be the set of moments corresponding to discrete spectral density functions. Define the set
\begin{equation}
\begin{split}
\Dscr:=\left\lbrace\,(\Psi,\Qb)\,: \Psi(\zetab_{\lb})>0,\ (\Psi^{-1}+Q)(\zetab_{\lb})>0,\ \right. \\  \left.\forall \lb \in \Zbb_\Nb^d\,\right\rbrace.
\end{split}
\end{equation}
Due to discretization, both $\Psi$ and $\Qb$ live in a finite dimensional vector space.
Consider the map
\begin{equation}\label{mmt_map}
\begin{split}
\fb:\, & \Dscr \to \Mscr_+ \\
& (\Psi,\Qb) \mapsto \Gamma ( (\Psi^{-1}+Q)^{-1} ).
\end{split}
\end{equation}
Given $\Sigmab\in\Mscr_+$, we aim to solve the equation
\begin{equation}\label{func_eqn}
\fb(\Psi,\Qb)=\Sigmab.
\end{equation}
When $\Psi$ is fixed, this is in fact equivalent to the stationarity condition $\nabla J_\Psi(\Qb)=0$ of the dual function (\ref{J_dual}), as explained in Subsection \ref{subsec:uniqueness}. Moreover, according to Theorem \ref{thm_main}, given a prior $\Psi$ and the moments $\{\Sigma_{\kb}\}_{\kb\in\Lambda}$ that are feasible, the solution $\Qb$ exists and is unique, and it can be obtained by minimizing the dual function (\ref{J_dual}). In other words, the solution map
\begin{equation*}
\sbf:\,(\Psi,\Sigmab) \mapsto \Qb
\end{equation*}
is well defined. 
We will next show that the solution map is smooth in either one of the two arguments when the other one is held fixed, as a consequence of Theorem \ref{thm_main}. The proof of the proposition below is deferred to the Appendix.

\begin{proposition}\label{prop_f_smooth}
	The map $\fb$ is smooth $($of class $C^\infty)$ on its domain $\Dscr$.
\end{proposition}


For a fixed $\Psi$ that is positive definite on $\Tbb_\Nb^d$, define the section of the map
\begin{equation}\label{omega_map}
\omega(\,\cdot\,):= \fb(\Psi,\,\cdot\,):\,\Lscr_+ \to \Mscr_+.
\end{equation}

\begin{theorem}\label{thm_cont_Sigma}
	The map $\omega$ is a diffeomorphism.
\end{theorem}

\begin{pf}
	Given a $\Sigmab\in\Mscr_+$, a solution to the equation $\omega(\Qb)=\Sigmab$ is a stationary point of the dual function $J_\Psi$ according to \eqref{stationary_cond}. By Theorem \ref{thm_main}, such a stationary point exists and is unique. Therefore, the map $\omega$ is a bijection. From the definition of $\omega$, one has $D\omega(\Qb)=D_2 \fb(\Psi,\Qb)$, i.e. the derivative of $\omega$ is equal to the partial derivative of $\fb$ with respect to the second argument. It then follows from Proposition \ref{prop_f_smooth} that $\omega$ is a smooth function.
	
	It now remains to prove the smoothness of the inverse $\omega^{-1}$, and this is an easy consequence of the inverse function theorem \citep[cf.][Theorem 5.2, p.~15]{lang1999fundamentals}. In order to see it, just notice that $D\omega(\Qb)$ is equal to the Hessian of the cost function $J_\Psi$ at $\Qb$ except for a sign difference. From Subsection \ref{subsec:uniqueness}, we know that the Hessian is positive definite for any $\Qb\in\Lscr_+$. Hence $D\omega(\Qb)$ is certainly a vector space isomorphism. By the inverse function theorem, we can conclude that $\omega$ is a local diffeomorphism at $\Qb$, and this implies the smoothness of $\omega^{-1}$.
\end{pf}

\begin{remark} \label{remark_1}
	As a consequence of Theorem \ref{thm_cont_Sigma} we have that for a fixed $\Psi$, the map $\sbf(\Psi,\,\cdot\,):\,\Mscr_+ \to\Lscr_+$ is continuous. Theorem \ref{thm_cont_Sigma} also implies that the set $\Mscr_+$ of moments indexed by $\Lambda$ is open. 
\end{remark}

We shall next show the well-posedness in the other argument, namely continuity of the map
\begin{equation}\label{imp_func_of_psi}
\sbf(\,\cdot\,,\Sigmab):\, \Psi \mapsto \Qb
\end{equation} 
when $\Sigmab$ is held fixed. Clearly, it is equivalent to solving the functional equation (\ref{func_eqn}) for $\Qb$ in terms of $\Psi$ when its right-hand side is fixed, which naturally falls in to the scope of the implicit function theorem.

\begin{theorem}\label{thm_cont_Psi}
	For a fixed $\Sigmab\in\Mscr_+$, the implicit function $\sbf(\,\cdot\,,\Sigmab)$ in $(\ref{imp_func_of_psi})$ is smooth.
\end{theorem}

\begin{pf}
	Fix a $\Psi$ that is positive definite on $\Tbb_\Nb^d$ and let $\Qb$ be the solution to \eqref{func_eqn}. Since the function $\Psi^{-1}+Q$ is positive definite on $\Tbb_\Nb^d$, it is not difficult to argue that there exist (open) neighborhoods $\mathcal U,\,\mathcal V$ of $\Psi$ and $\Qb$, respectively, such that $\tilde{\Psi}^{-1}+\tilde{Q}$ remains positive definite on $\Tbb_\Nb^d$ for any $\tilde{\Psi}\in \mathcal U$ and $\tilde{\Qb}\in \mathcal V$ using Lemma \ref{lem_bounded_below} in Appendix. We can therefore consider the function $\fb$ restricted to $\mathcal U\times \mathcal V$.
	
	The assertion then follows directly from the implicit function theorem \citep[see, e.g.,][Theorem 5.9, p.~19]{lang1999fundamentals}.
	More precisely, since the partial $D_2 \fb(\Psi,\Qb)$ is a vector space isomorphism, there exists a smooth map $\gb:\,\mathcal  U_0 \to \mathcal V$ defined on a sufficiently small open ball $\mathcal U_0 \subset \mathcal  U$, such that $\gb(\Psi)=\Qb$ and 
	\[ \fb( \tilde{\Psi} , \gb(\tilde{\Psi}) ) = \Sigmab \]
	for all $\tilde{\Psi}\in \mathcal  U_0$. Because there is a unique solution $\tilde{\Qb}$ corresponding to each $\tilde{\Psi}$, the restriction of $\sbf(\,\cdot\,,\Sigmab)$ on $\mathcal  U_0$ must coincide with $\gb$, and thus is smooth.
\end{pf}

In the next section we introduce a spectral estimation paradigm based on problem (\ref{primal_discrete}). Thanks to the results just proven, we
will establish the consistency of such an estimator.

\section{Spectral estimator}\label{sec:spec_est}
Let $\yb(\tb)$ be a second-order stationary random complex $m$-vector field. Next, we propose a spectral estimation procedure which uses a finite-size realization of the random field (i.e., dataset):
\begin{align}\label{dataset}
\mathcal Y=\{ y(\tb)\, ,\, 0\leq t_j \leq N_j-1 \hbox{ for } j=1,\dots,d\},
\end{align}
where $N_1,\dots,N_d$ are some positive integers. Since we do not have observations outside the index set $\Zbb^d_\Nb$ defined in \eqref{Zd_N}, we assume that the underlying random field is $\Nb$-periodic (cf.~Subsection~\ref{subsec:spec_rep}) whose spectral density is denoted by $\Phi(\zetab_{\lb})$ such that $\Phi(\zetab_{\lb})>0$ for all $\lb\in\Zbb^d_\Nb$.

In view of Section \ref{sec:problem}, a possible spectral estimate is obtained as follows:
\begin{enumerate}
	\item Set 
	\begin{equation}\label{Lambda_box}
	\Lambda=\{\kb\in\Zbb^d\,:\,|k_j|\leq n_j,\ j=1,\dots,d\},
	\end{equation}
	such that $N_j>2n_j$ with $j=1\ldots d$;
	\item Find a set of estimates $\{\hat \Sigma_\kb, \, \kb\in \Lambda\}$ of $\{ \Sigma_\kb, \, \kb\in \Lambda\}$ for which Assumption \ref{assump_feasibility} holds;
	\item An estimate of $\Phi(\zetab_{\lb})$ is the solution to (\ref{primal_discrete}) where $\Sigma_\kb$ has been replaced with $\hat \Sigma_\kb$. 
\end{enumerate}
Notice that $\Lambda$ is a $d$-dimensional box and condition $N_j>2n_j$ guarantees that Assumption \ref{assump_N} holds. The estimates $\{\hat \Sigma_\kb, \, \kb\in \Lambda\}$ may be computed from the data by taking the sample covariance and clearly they are ``reliable'' provided that $N_j\gg n_j$, i.e., the dataset $\mathcal Y$ is long enough along each dimension.  
In what follows, the spectral estimator obtained as above will be denoted by $\hat \Phi(e^{i\thetab})$.

The remaining nontrivial step is to construct the set of estimates $\{\hat \Sigma_\kb, \, \kb\in \Lambda\}$ satisfying Assumption \ref{assump_feasibility}. The latter
plays an important role in the previous proofs concerning the well-posedness of the optimization problem \eqref{primal_discrete}. In the unidimensional case ($d=1$) with the normalized Lebesgue measure, feasibility can be checked and enforced \citep{Zorzi-F-12,FPZ-12} as it is equivalent to a simple algebraic condition. More precisely, given a finite sequence of estimates $\hat\Sigma_0,\hat\Sigma_1,\dots,\hat\Sigma_n$, the set of moment equations
\begin{equation}
\int_{\Tbb} e^{ik\theta}\Phi(e^{i\theta})\frac{\d\theta}{2\pi}=\hat\Sigma_k,\quad k=0,\dots,n
\end{equation}
has a solution $\Phi:\,\Tbb\to\Hfrak_{+,m}$ if and only if the block-Toeplitz matrix
\begin{equation}
\bmat \hat\Sigma_0&\hat\Sigma_1^*&\hat\Sigma_2^*&\cdots&\hat\Sigma_n^* \\
\hat\Sigma_1&\hat\Sigma_0&\hat\Sigma_1^*&\cdots&\hat\Sigma_{n-1}^* \\
\hat\Sigma_2&\hat\Sigma_1&\hat\Sigma_0&\cdots&\hat\Sigma_{n-2}^* \\
\vdots&\vdots&\ddots&\ddots&\vdots \\
\hat\Sigma_n&\hat\Sigma_{n-1}&\cdots&\hat\Sigma_1&\hat\Sigma_0\emat
\end{equation}
is positive definite. The proof of this fact is intrinsically related to the factorization problem of positive matrix Laurent polynomials. A similar positivity condition exists in high dimensional cases. For example, in the $2$-d case, a block-Toeplitz matrix with Toeplitz blocks is required to be positive definite \citep[cf.][]{geronimo2004positive}.
However, when the dimension $d>1$, such a positivity condition is only a necessary condition for the existence of a solution to the moment equations, and a sufficient condition seems yet unknown, partly due to the fact that the factorization problem is quite difficult even for scalar polynomials in several variables \citep[see e.g.,][]{geronimo2006factorization}.


In view of the difficulty of attacking the feasibility question directly, in this section we will instead give a method that guarantees the feasibility of our optimization problem. Essentially, this is a multivariate and multidimensional generalization of the \emph{standard biased covariance estimates} for scalar unidimensional processes reported in \citet[Chapter~2]{stoica2005spectral}. Define the finite Fourier transform
\begin{equation}\label{y_sample_FT}
\hat{\yb}_{\Nb}(e^{i\thetab}) := \sum_{\tb\in\Zbb^d_{\Nb}} \yb(\tb) e^{-i\innerprod{\tb}{\thetab}}.
\end{equation}
Then the \emph{periodogram} is defined as
\begin{equation}\label{Phi_periodgram}
\hat{\Phi}_{\p}(e^{i\thetab}) := \frac{1}{|\Nb|} \hat{\yb}_{\Nb}(e^{i\thetab}) \hat{\yb}_{\Nb}(e^{i\thetab})^*+\frac{\varepsilon}{|\Nb|}I
\end{equation}
where
the real constant $\varepsilon>0$ can be arbitrarily small. It is worth noting that the first term on the right hand side of (\ref{Phi_periodgram}) is rank one and positive semidefinite. Therefore, the bias term $\varepsilon/|\Nb| I$ guarantees that $\hat{\Phi}_{\p}(e^{i\thetab})$ is positive definite. It is worth noting that $\hat{\Phi}_{\p}(e^{i\thetab})$ is a first estimate of $\Phi(e^{i\thetab})$ due to the relation 
\begin{equation}\label{Phi_def2}
\Phi(e^{i\thetab}) = \lim_{\min(\Nb)\to\infty} \E \hat{\Phi}_{\p}(e^{i\thetab})
\end{equation}
which holds under a mild assumption on the decay rate of the covariance lags. The more precise statement is given in the next proposition. To this aim we  introduce the set
\begin{equation}\label{Zd_2K}
\Zbb^d_{2\Nb-\oneb} := \{ \lb\in\Zbb^d \,:\, -N_j+1 \le \ell_j \le N_j-1,\, j=1,\dots,d \}
\end{equation}
for the covariance lags.

\begin{proposition}\label{prop_equiv_def_psd}
	The equality \eqref{Phi_def2} holds if
	\begin{equation}\label{cond_cov_decay}
	\lim_{\min(\Nb)\to\infty} \frac{1}{|\Nb|} \sum_{\kb\in\Zbb^d_{2\Nb-\oneb}} \|\Sigma_{\kb}\| r_{\kb} = 0,
	\end{equation}
	where $r_{\kb}:=|\Nb|-\prod_{j=1}^d \left( N_j-|k_j| \right)$.
\end{proposition}

\begin{pf}
	Given the expression \eqref{y_sample_FT}, it is easy to compute
	\begin{equation}\label{exp_Phi_p}
	\begin{split}
	\E \hat{\Phi}_{\p}(e^{i\thetab}) & = \frac{1}{|\Nb|} \sum_{\tb\in\Zbb^d_{\Nb}} \sum_{\sbf\in\Zbb^d_{\Nb}} \Sigma_{\tb-\sbf} e^{-i\innerprod{\tb-\sbf}{\thetab}} +\frac{\varepsilon}{|\Nb|}I\\
	& = \frac{1}{|\Nb|} \sum_{\kb\in\Zbb^d_{2\Nb-\oneb}} \Sigma_{\kb} e^{-i\innerprod{\kb}{\thetab}} c_{\kb}+\frac{\varepsilon}{|\Nb|}I,
	\end{split}
	\end{equation}
	where the index $\kb=\tb-\sbf$ results from the substitution, and $\{c_{\kb}\}$ are some integers to be determined. The index pair $(\tb,\sbf)$ contributing to $c_{\kb}$ must satisfy the inequalities
	\begin{equation}
	\begin{split}
	0 \le & s_j \le N_j-1 \\
	0 \le & t_j=s_j+k_j \le N_j-1
	\end{split}
	\end{equation}
	for $j=1,\dots,d$, which in turn yield
	\begin{equation}\label{s_j_inequal}
	\left\{ \begin{array}{ll}
	0 \le s_j \le N_j-1-k_j & \textrm{if } k_j \ge 0,\\
	-k_j \le s_j \le N_j-1 & \textrm{if } k_j<0.\\
	\end{array} \right.
	\end{equation}
	Therefore, the number of admissible $s_j$ is $N_j-|k_j|$, and consequently $c_\kb=\prod_{j=1}^d \left( N_j-|k_j| \right) = |\Nb|-r_\kb$. Now following \eqref{exp_Phi_p}, we have
	\begin{align}
	\lim_{\min(\Nb)\to\infty} \E \hat{\Phi}_{\p}(e^{i\thetab}) = \lim_{\min(\Nb)\to\infty} \left\{ \sum_{\kb\in\Zbb^d_{2\Nb-\oneb}} \Sigma_{\kb} e^{-i\innerprod{\kb}{\thetab}}\right.\nonumber \\  \left. - \frac{1}{|\Nb|} \sum_{\kb\in\Zbb^d_{2\Nb-\oneb}} \Sigma_{\kb} e^{-i\innerprod{\kb}{\thetab}} r_{\kb} +\frac{\varepsilon}{|\Nb|}I\right\}\nonumber 
	\end{align}
	which is equal to $\Phi(e^{i\thetab})$ defined in \eqref{Phi_spec_density} under the stated assumption \eqref{cond_cov_decay}.
\end{pf}

It is worth noting that $\hat \Phi_\p(\zetab_{\lb})$ is a crude estimator of $\Phi(\zetab_{\lb})$ because it is not consistent; indeed, from (\ref{Phi_periodgram}) we have that if $\hat\Phi_\p(\zetab_{\lb})$ converges to a deterministic quantity as $\min(\Nb)\rightarrow \infty$ then the latter is a rank one positive semidefinite matrix. In contrast, we have $\Phi(\zetab_{\lb})>0$ by assumption.   

Next we shall derive those covariance estimates that correspond to the periodogram. We have:
\begin{equation} 
\begin{split}
\hat{\Phi}_{\p}&(e^{i\thetab})  = \frac{1}{|\Nb|} \sum_{\tb\in\Zbb^d_{\Nb}} \sum_{\sbf\in\Zbb^d_{\Nb}} \yb(\tb) \yb(\sbf)^* e^{-i\innerprod{\tb-\sbf}{\thetab}}+\frac{\varepsilon}{|\Nb|}I \\
& = \frac{1}{|\Nb|} \sum_{\kb\in\Zbb^d_{2\Nb-\oneb}} \sum_{\sbf\in\Xi_{\Nb,\kb}} \yb(\sbf+\kb) \yb(\sbf)^* e^{-i\innerprod{\kb}{\thetab}}+\frac{\varepsilon}{|\Nb|}I,
\end{split}
\end{equation}
where the index set
\begin{equation}
\Xi_{\Nb,\kb}:=\{\sbf\in\Zbb^d \,:\, s_j \textrm{ satisfies \eqref{s_j_inequal} for } j=1,\dots,d\}.
\end{equation}
For $\kb\in\Lambda$, an estimator of $\Sigma_\kb$ is given by
\begin{align} \label{estim_lag}
\hat  \Sigma_{\kb}= \int_{\Tbb^d}e^{i\innerprod{\kb}{\thetab}}\hat \Phi_\p(e^{i\thetab})\d\nu_\Nb(\thetab) .
\end{align}
Clearly, such a set of estimates satisfies Assumption \ref{assump_feasibility}.
We conclude this section by showing that $\hat \Phi(e^{i\thetab})$ is consistent.
\begin{proposition} Consider the parametric family of spectral densities:
	\begin{align}
	\mathcal F_\Psi=\{(\Psi^{-1}+Q)^{-1} \hbox{ s.t. } (\Psi^{-1}+Q)(\zetab_{\lb})>0 \, \forall \lb \in \Zbb_\Nb^d\}\nonumber
	\end{align}
	where $Q(e^{i\thetab})=\sum_{\kb\in\Lambda} Q_\kb e^{-i<\kb,\thetab>}$, with $Q_\kb\in\mathbb C^{m\times m}$, and $\Psi(\zetab_{\lb})>0$ $\forall \lb \in \Zbb_\Nb^d$ is fixed. Assume that the dataset \eqref{dataset} is generated by an ergodic random field having spectral density $\Phi(e^{i\thetab})\in\mathcal F_\Psi$. Let $\hat \Phi(e^{i\thetab})$  denote the IS spectral estimator obtained using the same prior $\Psi$ and the covariance estimates $(\ref{estim_lag})$, then we have 
	\begin{align}\label{claim_er}
	\hat \Phi(e^{i\thetab})\overset{\hbox{\footnotesize a.s.}}{\longrightarrow}\Phi(e^{i\thetab})
	\end{align}
	as $\min(\Nb)\rightarrow \infty$.
\end{proposition}

\begin{pf}
	Since $\Phi(e^{i\thetab})\in\mathcal F_\Psi$, there exists $\tilde Q(e^{i\thetab})=\sum_{\kb\in\Lambda} \tilde Q_\kb e^{-i<\kb,\thetab>}$ such that $\Phi=(\Psi+\tilde Q)^{-1}$. In particular, we have 
	\begin{align}\label{prop_er1}
	\sbf(\Psi,\Sigmab)=\tilde Q
	\end{align}
	with $\Sigmab=\Gamma(\Phi)$. Since the random field is ergodic we have 
	\begin{align}\label{prop_er2}
	\hat \Sigma_\kb\overset{\hbox{\footnotesize a.s.}}{\longrightarrow} \Sigma_\kb, \; \forall \, \kb\in\Lambda
	\end{align} as $\min(\Nb)\rightarrow \infty$.
	Let $Q_{(\Nb)}$ be the solution to the dual problem (\ref{J_dual}) where we have made explicit its dependence with respect to the size of the dataset $\mathcal Y$. Then, $\hat \Phi=(\Psi^{-1}+Q_{(\Nb)})^{-1}$ and \begin{align}
	\sbf(\Psi,\hat \Sigmab)=Q_{(\Nb)}
	\end{align} where $\hat \Sigmab=\Gamma(\hat\Phi_\p)$. Since $\sbf$ is continuous with respect to the second argument (see Remark \ref{remark_1}) and in view of (\ref{prop_er1})--(\ref{prop_er2}) we have that 
	\begin{align}
	Q_{(\Nb)} \overset{\hbox{\footnotesize a.s.}}{\longrightarrow}\tilde Q 
	\end{align} and thus we obtain (\ref{claim_er}).
\end{pf}

\section{Numerical examples}\label{sec:simulation}
In this section, we apply our theory to the problem of target parameter estimation in an integrated system of two automotive modules, see Figure \ref{fig_auto}. 
\begin{figure}[!t]
	\centering
	\includegraphics[width=0.48\textwidth]{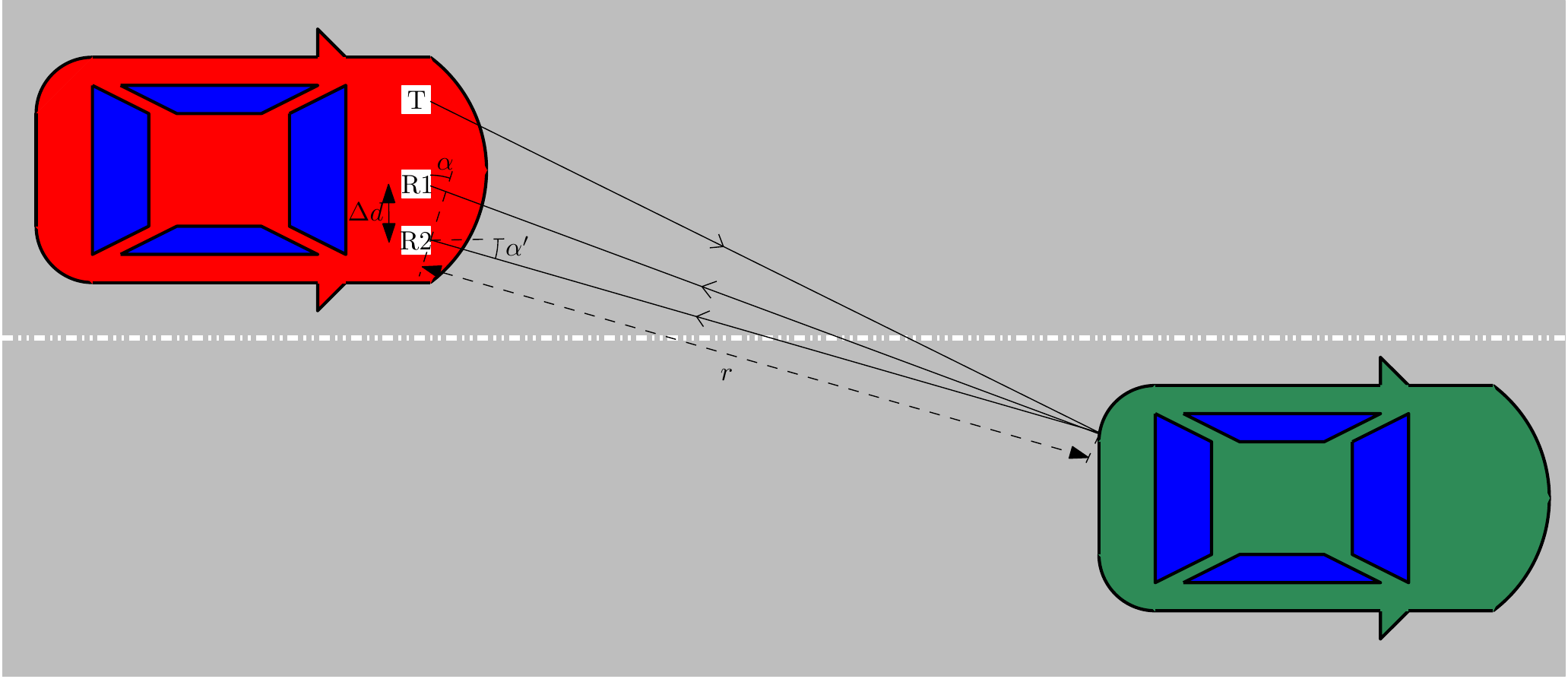}
	\caption{Integrated system of two automotive modules: T is the transmitter, R1 and R2 are two ULAs (receivers).}
	\label{fig_auto}
\end{figure}
The setting of the problem is the same as that in \citet{ZFKZ2019fusion} which we will recall briefly. For details of the radar signal (the waveform, filtering, sampling, etc.), we refer readers to the literature \citet{rohling2012continuous,engels2014target,engels2017advances}. For simplicity, let us also assume that only one target is present in the field of view. Our measurements come from two uniform linear arrays (ULAs) of receive antennas denoted by R1 and R2, respectively, in Fig. \ref{fig_auto}. The two ULAs are placed in the same line at a distance $\Delta d$, say a few decimeters. In a very short time interval at the scale of 20 ms (called ``coherent processing interval''), the scalar measurement of each ULA is modeled as a $3$-d complex sinusoid \citep[cf.][]{engels2014target}, that is, for each vector $\tb\in\Zbb_\Nb^3$,
\begin{equation}\label{y_j_measurement}
\begin{split}
y_1(\tb) & = a \, e^{ i \left( \innerprod{\thetab}{\tb} +\varphi \right)} + w_1(\tb) \\
y_2(\tb) & = a \, e^{ i \left( \innerprod{\thetab}{\tb} + M\theta_3 +\varphi \right)} + w_2(\tb)
\end{split}
\end{equation}
where the subscripts $1,2$ label the ULAs. Model (\ref{y_j_measurement}) holds under the \emph{far field} assumption which is common in this kind of radar applications. The meaning of each variable will be explained next.
The index $\tb$ takes values in the set \eqref{Zd_N} with the dimension $d=3$, and the vector $\Nb=[N_1,N_2,N_3]$ corresponds to the size of the measurement. Here the integer $N_1$ denotes the number of samples per pulse, $N_2$ the number of pulses, and $N_3$ the number of (receive) antennas. The scalar $a$ is a real amplitude. 
The variable $\varphi$ is an initial phase angle of the first measurement channel which is assumed to be a random variable uniformly distributed in $[-\pi,\pi]$ \citep[cf.][Section~4.1]{stoica2005spectral}.
The processes $w_k$ $(k=1,2)$ are uncorrelated zero-mean circular complex white noises with the same variance ${\sigma}^2$, and both are independent of $\varphi$.
The real vector $\thetab=[\theta_1\; \theta_2\;\theta_3]\in\Tbb^3$ contains three \emph{unknown} normalized angular target frequencies. The components $\theta_j$ $(j=1,2,3)$ are related to the range $r$, the (radial) relative velocity $v$, and the azimuth angle $\alpha$ of the target via
simple invertible functions \citep[see][Section~16.4]{engels2014target}, such that the target parameter vector $(r,v,\alpha)$ can be readily recovered from the frequency vector $\thetab$.
The number $M=\Delta d/\Delta s$ where $\Delta s$ is the distance between two adjacent antennas in the ULA, and $M\theta_3$ represents the phase shift between the measurements of the two ULAs due to the distance $\Delta d$. The target parameter estimation problem consists in estimating the unknown target frequencies $\thetab$ from the sinusoid-in-noise measurements generated according to model \eqref{y_j_measurement}.

Such a frequency estimation problem has been extensively studied in the literature \citep[see, e.g.,][Chapter~4]{stoica2005spectral}, and many methods have been proposed to solve it, in the case of a single ULA.
We will address the problem via M$^2$ spectral estimation as explained next. 
{Notice that the current problem setup falls into our M$^2$ framework because two ULAs produce a multivariable (bivariate) signal and the three physical quantities of the target (azimuth angle, velocity, range) relating to the angular frequencies give rise to a three dimensional domain.}
Set $\yb(\tb):=[\,y_1(\tb)\,y_2(\tb)\,]^\top$. Through elementary calculations, we have
\begin{equation}
\Sigma_{\kb}:= \E \yb(\tb+\kb) \yb(\tb)^* = a^2 e^{i\innerprod{\thetab}{\kb}} R + {\sigma}^2 \delta_{\kb,\zerob} I_2,
\end{equation}
where the matrix
\begin{equation}\label{mat_R}
R=\bmat
1&e^{-iM\theta_3}\\
e^{iM\theta_3}&1
\emat,
\end{equation}
and $\delta_{\kb,\zerob}$ is the Kronecker delta function. Taking the Fourier transform, the multidimensional-multivariate spectrum of $\yb$ is
\begin{equation}\label{Phi_ideal}
\Phi(e^{i\omegab}) = 2\pi a^2 \delta(\omegab-\thetab) R + {\sigma}^2 I_2
\end{equation}
where $\delta(\cdot)$ here is the Dirac delta measure. Although the above spectrum is singular, the idea is to approximate it with a nonsingular spectrum with a peak in $\thetab$. Therefore, we first compute an estimate $\hat{\Phi}$ of the ideal spectrum from the radar measurements following the procedure described in Section \ref{sec:spec_est}. Then we use the post-processing method proposed in \citet{ZFKZ2019fusion} to obtain an estimate of the target frequency vector via
\begin{equation}\label{theta_est_m2}
\hat{\thetab}_\F:= \underset{\omegab\in\Tbb^3}{\argmax} \, \|\hat{\Phi}(e^{i\omegab})\|_{\F}^2,
\end{equation}
where the subscript $_\F$ of the estimate $\hat{\thetab}$ stands for ``Frobenius'' as
$\|\hat{\Phi}(e^{i\omegab})\|_{\F}^2:=|\hat{\Phi}_{11}(e^{i\omegab})|^2+|\hat{\Phi}_{22}(e^{i\omegab})|^2+2\,|\hat{\Phi}_{12}(e^{i\omegab})|^2$ is the Frobenius norm (squared). As discussed in \citet{ZFKZ2019fusion} the cross spectrum $\hat \Phi_{12}$ merges the information coming from the two measurement channels and improves the estimation of $\thetab$.   


We report below one numerical example in which the data size $\Nb=[30\;30\;8]$, the amplitude of the sinusoid $a=1$, the number $M=20$, and the noise variance ${\sigma}^2=2$. The grid size of the discrete $3$-torus is the same as $\Nb$. The true frequency vector for the data generation is $\thetab=[0.8101\;-0.5872\; 2.1798]$, while its quantized version on $\Tbb_\Nb^3$ is ${\thetab}_q=[0.8378\;-0.6283\; 2.3562]$ corresponding to the $3$-d index $[5\;28\;4]$.\footnote{The array index starts from $1$ under the convention of Matlab.} In what follows we consider the estimator of Section \ref{sec:spec_est}, denoted by IS, where the set $\Lambda$ is defined in \eqref{Lambda_box} with $\nb=[1\; 1\; 1]$. The prior is taken as $\Psi\equiv\hat{\Sigma}_\zerob$, the (constant) estimated zeroth moment. {The IS solver is initialized with $\Qb=\zerob$ which is feasible since we have $\hat{\Sigma}_\zerob>0$ according to the estimation scheme \eqref{estim_lag}.}
We compare the IS method with the two windowed M$^2$ periodograms proposed in \citet{ZFKZ2019fusion}, one with a rectangular window of size $[8\; 8\; 2]$, denoted by RECT,  and the other with a Bartlett window of size $[12\; 12\;  3]$, denoted by BART. Using the post-processing method in \eqref{theta_est_m2}, both the periodograms and our method return the same frequency estimate ${\thetab}_q$, which is the best grid point that approximates the true $\thetab$. However, from Figs. \ref{fig:res_main54_sec1}, \ref{fig:res_main54_sec2}, and \ref{fig:res_main54_sec3}, showing  the three sections of the function $\|\hat{\Phi}(e^{i\omegab})\|_{\F}^2$ for the different estimated spectra, only the IS estimator exhibits a proper peak {that is clearly distinguishable form possible background noise. The latter property is desirable for the peak detection task.}
\begin{figure}[!t]
	\centering
	\includegraphics[width=0.5\textwidth]{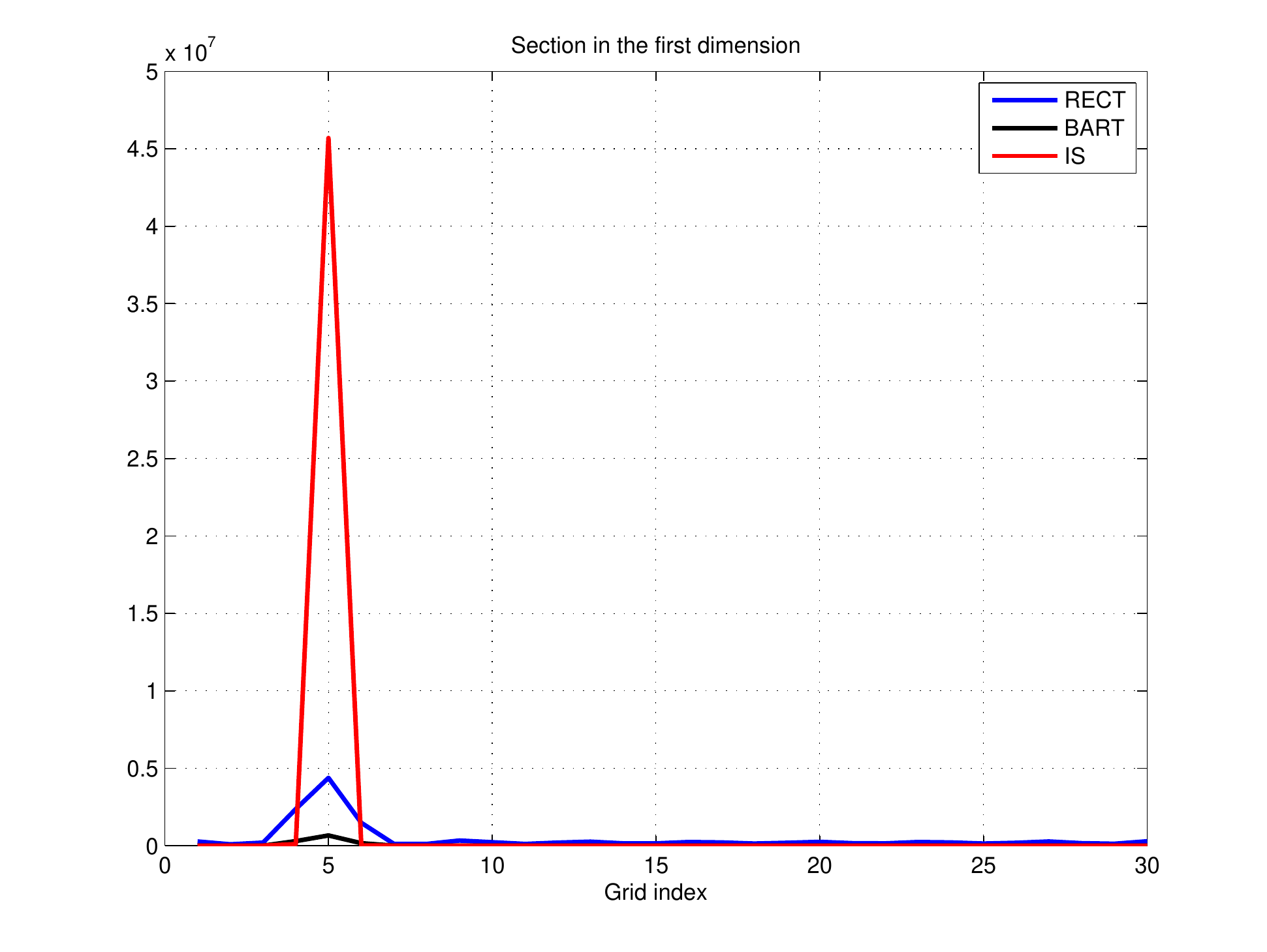}
	\caption{{Squared \emph{Frobenius} norms  of the estimated spectra at the cross section $[\,\cdot \;28\;4\,]$, i.e., $\|\hat{\Phi}(e^{i\omegab})\|_{\F}^2$ with $\omegab=2\pi\times[(k-1)/30, 27/30, 3/8]$ for grid indices $k=1,\ldots, 30$, where the IS estimator exhibits a proper peak.}}
	\label{fig:res_main54_sec1}
\end{figure}
\begin{figure}[!t]
	\centering
	\includegraphics[width=0.5\textwidth]{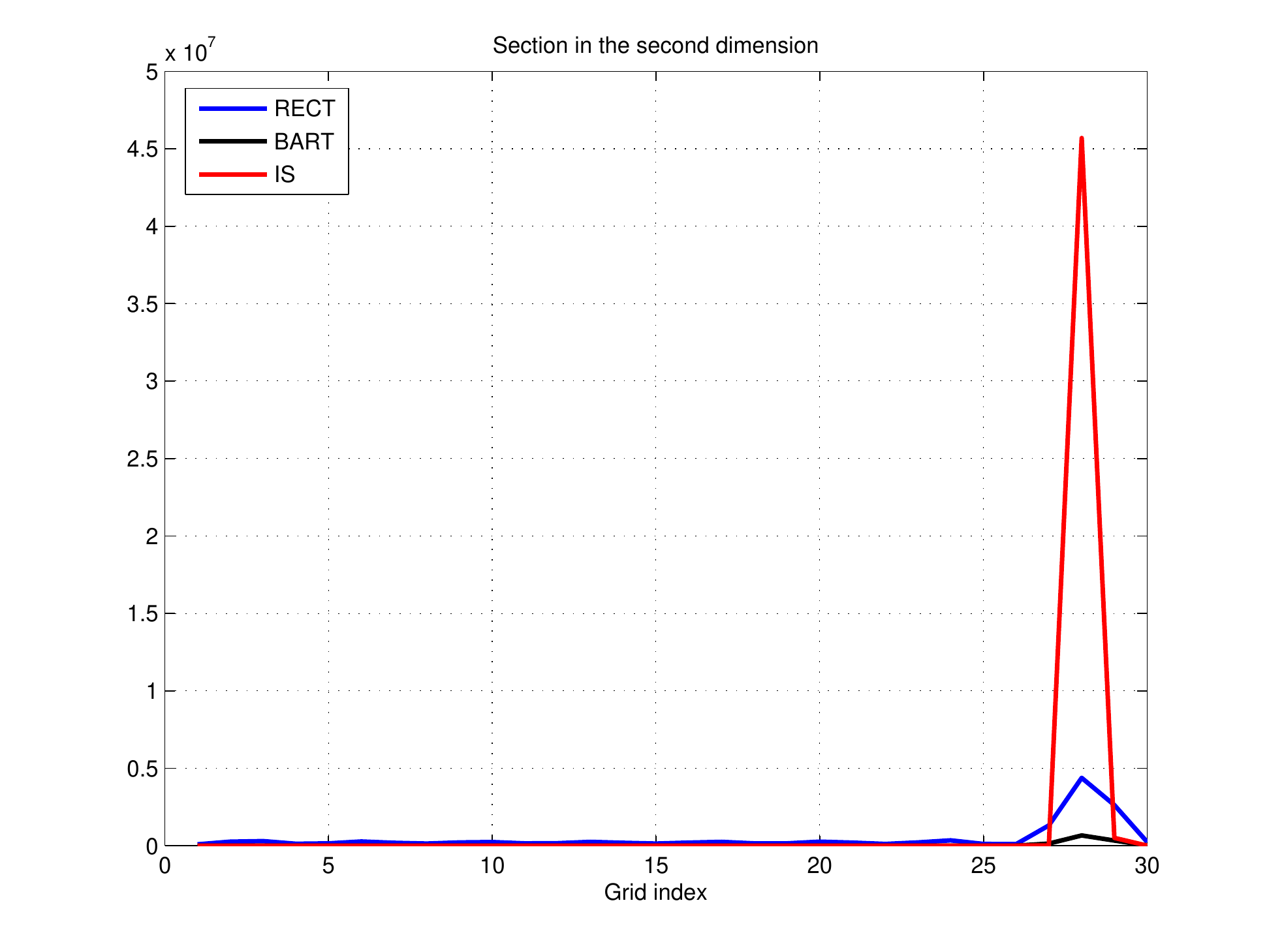}
	\caption{{Squared \emph{Frobenius} norms  of the estimated spectra at the cross section $[\, 5 \;\cdot\;4\,]$, i.e., $\|\hat{\Phi}(e^{i\omegab})\|_{\F}^2$ with $\omegab=2\pi\times[4/30, (k-1)/30, 3/8]$ for grid indices $k=1,\ldots, 30$, where the IS estimator exhibits a proper peak.}}
	\label{fig:res_main54_sec2}
\end{figure}
\begin{figure}[!t]
	\centering
	\includegraphics[width=0.5\textwidth]{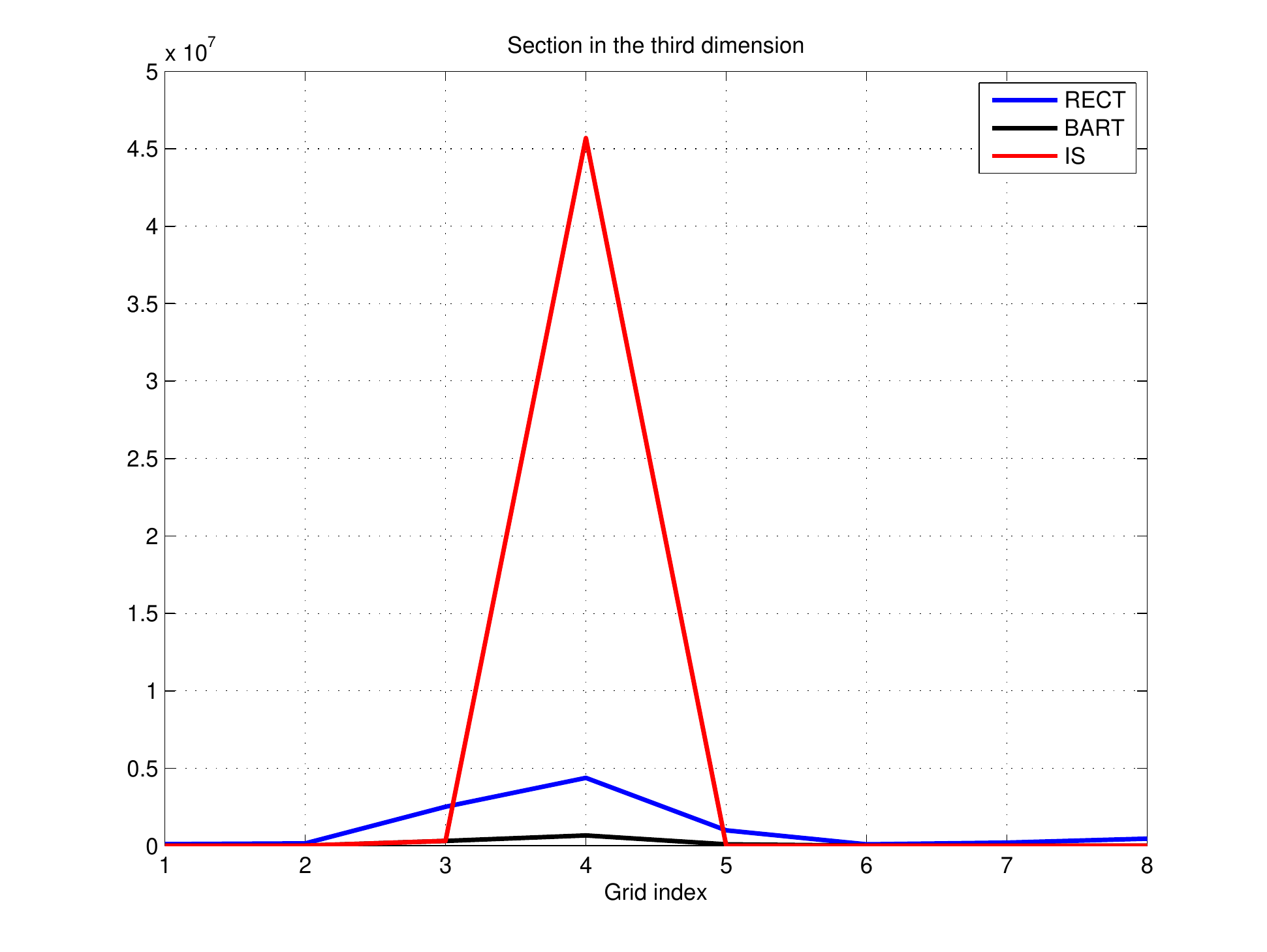}
	\caption{{Squared \emph{Frobenius} norms  of the estimated spectra at the cross section $[\,5\; 28\;\cdot\,]$, i.e., $\|\hat{\Phi}(e^{i\omegab})\|_{\F}^2$ with $\omegab=2\pi\times[4/30, 27/30, (k-1)/8]$ for grid indices $k=1,\ldots, 8$, where the IS estimator exhibits a proper peak.}}
	\label{fig:res_main54_sec3}
\end{figure} We have also run Monte-Carlo simulations for this frequency estimation task. Under the parameter configuration above, the IS method performs in the same way as the M$^2$ periodogram-based spectral estimators as measured by the error $\|\hat{\thetab}-\thetab\|$.

Since the true spectrum in the above radar application contains only spectral lines, it does not belong to the model class which the IS method produces. Our method gives a rational approximation of the spectral line, but it is difficult to quantify how good such an approximation is due to the singularity of the true spectrum. Out of such consideration, we want to test our method when the generative model is rational. Consider the autoregressive (AR) model 
\begin{equation}\label{AR_uni}
x(t) = \frac{1}{1-\alpha z^{-1}} w(t),
\end{equation}
where $w(t)$ is a white noise with unit variance and $\alpha=\rho\, e^{i\theta}$ is such that the modulus $\rho$ is close to $1$ and $\theta\in[-\pi,\pi]$.
It is well-known that in the scalar unidimensional case, the above AR model approximates the sinusoid in the sense that the AR spectrum has a peak at frequency $\theta$. A possible generalization of (\ref{AR_uni}) in the multidimensional case is 
\begin{equation}\label{x_AR}
x(\tb) = \frac{1}{1-\innerprod{\alphab}{\zb^{-\oneb}}} w(\tb)
\end{equation}
where $\alpha_j=\rho_j \, e^{i\theta_j}\ (j=1,2,3)$ and $\innerprod{\alphab}{\zb^{-\oneb}}:=\alpha_1 z_1^{-1}+\alpha_2 z_2^{-1}+\alpha_3 z_3^{-1}$ are such that the sum $\rho_1+\rho_2+\rho_3$ is close to $1$. The peak of the spectrum is of course obtained at the vector $\thetab:=[\theta_1,\theta_2,\theta_3]$ of phase angles. We give next the result of a Monte-Carlo simulation that contains $100$ trials. In each trial, each component of the frequency vector $\thetab$ is generated from the uniform distribution in $[-\pi,\pi]$. The signal model for the measurement is
\begin{equation}\label{y_j_ARinNoise}
\begin{split}
y_1(\tb) & = x(\tb) + w_1(\tb) \\
y_2(\tb) & = x(\tb) \, e^{ i M\theta_3} + w_2(\tb)
\end{split}
\end{equation}
which mimics the sinusoid-in-noise model \eqref{y_j_measurement}. The difference here is that the true signal is replaced with the AR process defined by \eqref{x_AR}, in which we have chosen the pole moduli $\rho_1=\rho_2=\rho_3=0.3$. The variance ${\sigma}^2$ of the additive noise in \eqref{y_j_ARinNoise} is $2$. 
{The realization of the process $x(\tb)$ is generated by applying the $3$-d recursion \eqref{x_AR} given the noise and the zero boundary condition. Notice that in order to reach the ``steady state'', a realization of size $1000^3$ is generated recursively and only the last $\Nb$ samples are retained for the covariance estimation.}
The true spectrum $\Phi$ of the process \eqref{y_j_ARinNoise} is given by
\begin{equation}\label{Phi_true}
\Phi(e^{i\omegab}) = \Phi_x(e^{i\omegab}) R + {\sigma}^2 I_2
\end{equation}
which is certainly rational with
\begin{equation}
\Phi_x (e^{i\omegab}) = \frac{\mathrm{Var}(w)}{\left|1-\innerprod{\alphab}{e^{-i\omegab}}\right|^2} ,
\end{equation}
the spectral density of the AR process \eqref{x_AR}.
We are interested in how well the true spectrum $\Phi$ is approximated by the IS estimator and the M$^2$ periodograms. Let us define the relative error $\|\hat{\Phi}-\Phi\| / \|\Phi\|$ where $\hat{\Phi}$ is one of the spectrum estimates. These errors of the different methods in $100$ trials are reported in Fig.~\ref{fig:res_main83_apx}. The data size $\Nb$, the discrete torus $\Tbb^3_\Nb$, the model order $\nb$ and the prior $\Psi$ of the IS method, as well as the window widths of the periodograms are the same as those in the previous part concerning sinusoidal signals.
One can see that the IS estimator clearly outperforms the periodograms.  
\begin{figure}[!t]
	\centering
	\includegraphics[width=0.5\textwidth]{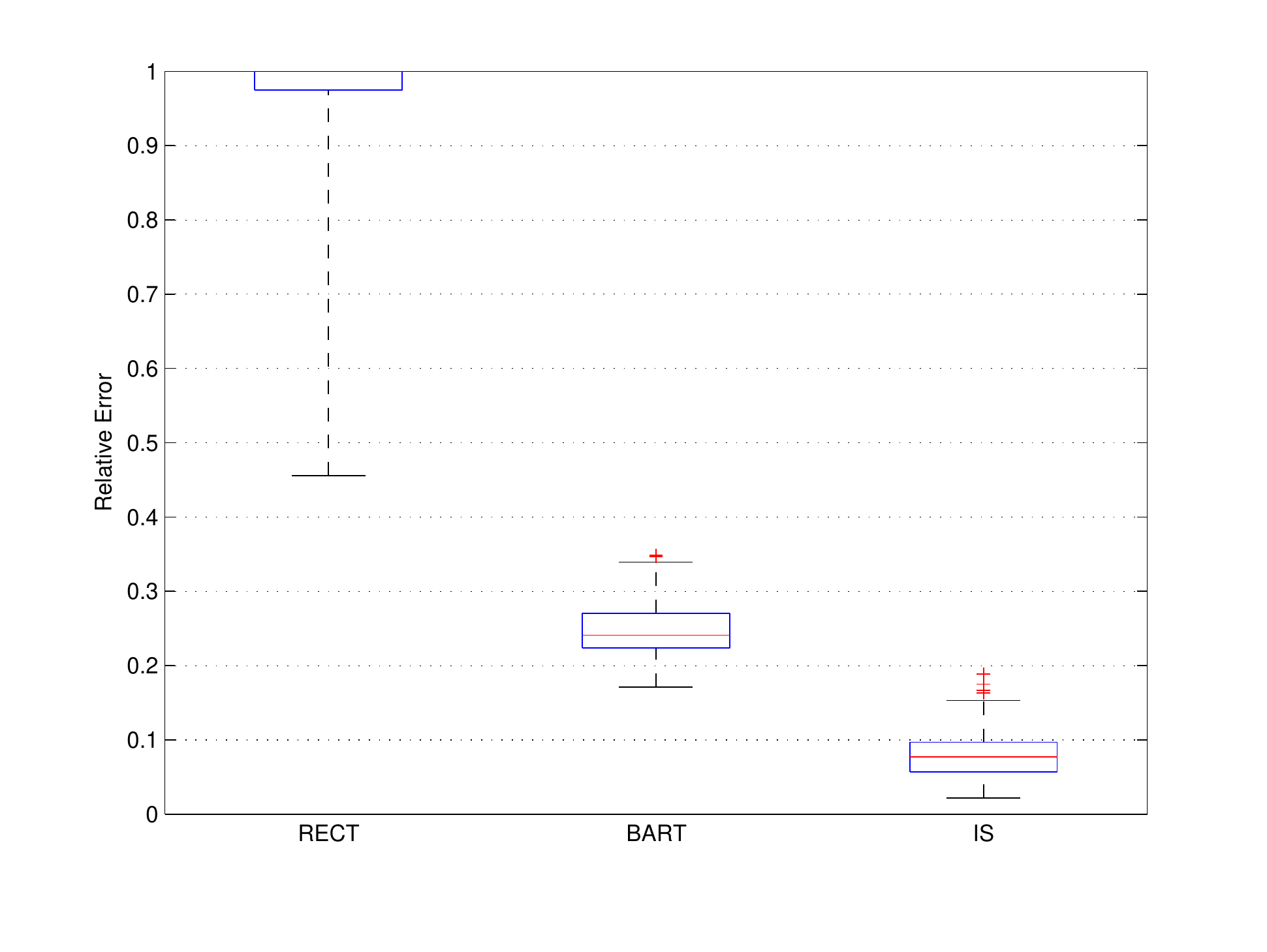}
	\caption{Relative error of spectrum approximation using different methods. Note that it is not possible show the entire boxplot of RECT because its performance is very poor compared to BART and IS.}
	\label{fig:res_main83_apx}
\end{figure}
Moreover, we also want to compare different methods of spectral estimation in peak finding, namely how well the estimators can locate the peak of the true spectrum at $\thetab$.
The errors $\|\hat{\thetab}-\thetab\|$ of the estimated peak location $\hat{\thetab}$ returned by different methods during the $100$ trials are depicted in Fig.~\ref{fig:res_main83_est}. The three boxplots indicate that the IS method also outperforms the M$^2$ periodograms in peak finding, which is within our expectation given the result in Fig.~\ref{fig:res_main83_apx}. {As for the computational speed, the current implementation to solve the optimization problem associated with the IS estimator is slower than the periodograms since the latter involves essentially only (linear) FFT operations.}
\begin{figure}[!t]
	\centering
	\includegraphics[width=0.5\textwidth]{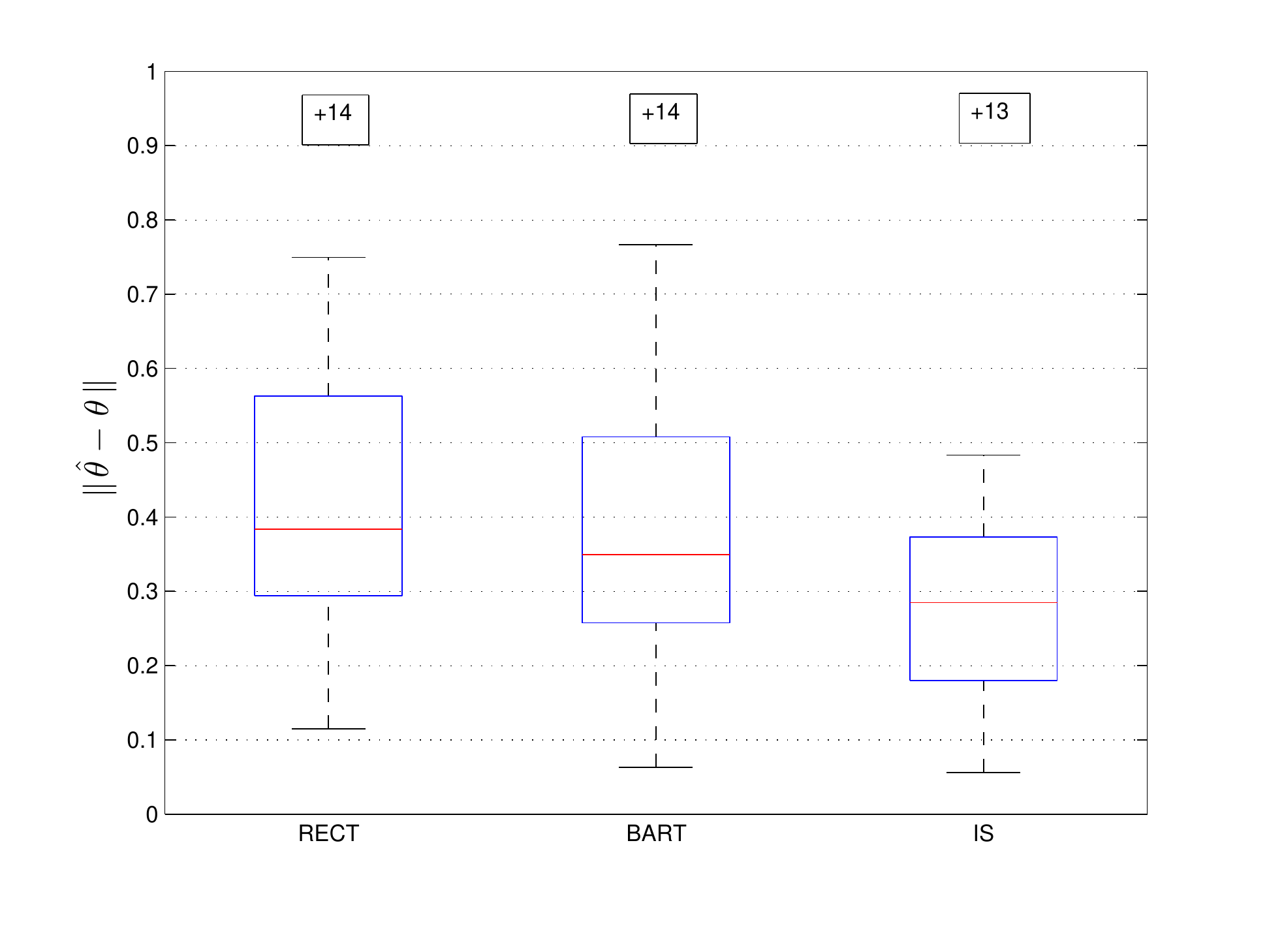}
	\caption{Error of peak finding using different methods.  The values on the boxes represent the number of outliers greater than 1 for each estimator.}
	\label{fig:res_main83_est}
\end{figure}

\section{Conclusions}
\label{sec:conclusions}

In this paper, we have considered a M$^2$ spectral estimation problem for periodic stationary random vector fields via covariance extension, i.e., matching a finite number of matrix covariance lags. Our solution is based on convex optimization with the Itakura-Saito distance which incorporates an \emph{a priori} spectral density as extra information. We have shown that the optimization problem is well-posed, and thus a smooth parametrization of solutions can be obtained by changing the prior function. Moreover, a covariance estimation scheme has been proposed given a finite-size realization of the random field, and a spectrum estimation procedure using these estimated covariances is described. To illustrate the theory, we have performed numerical simulations concerning the parameter estimation problem in an automotive radar system. The results show that our spectral estimator is very competitive with periodogram-based spectral estimators. 
{For practical applications however, efficient methods need to be developed so that one can compute IS estimators in real time, and this will be the subject of further study.
As another future research topic}, we plan to attack a spectral estimation problem similar to \eqref{primal_cont} in which the spectrum is define on the continuum of $\Tbb^d$. Such a problem will no doubt be much more challenging, and careful analysis must be carried out on the boundary of the feasible set.


\bibliographystyle{model4-names}
\bibliography{references}           

\appendix

\section{Proof of Proposition~\ref{prop_f_smooth}}

The next lemma is need in the proof of Proposition \ref{prop_f_smooth}.

\begin{lemma}\label{lem_bounded_below}
	Let $M$ be a function on $\Tbb_\Nb^d$ whose values are Hermitian positive definite matrices. If another function $\tilde{M}$ is sufficiently close to $M$ in norm, then there exists a real constant $\mu>0$ such that $\tilde{M}(\zetab_{\lb}) \geq \mu I$ for all $\lb\in\Zbb_\Nb^d$.
\end{lemma}
\begin{pf}
	Due to the fact that we are considering functions on the grid $\Tbb_\Nb^d$, the lemma follows directly from the continuous dependence of the eigenvalues on the matrix \citep[cf.][Theorem 6, p.~130]{lax2007linear}.
\end{pf}

\begin{pf}[Proof of Proposition $\ref{prop_f_smooth}$]
	Let us first show that $\fb$ is of class $C^1$. According to \citet[Propositions 3.4 \& 3.5, p.~10]{lang1999fundamentals}, it is equivalent to show that the two partial derivatives of each ``component''
	\begin{equation}\label{f_kb}
	\fb_\kb:\, (\Psi,\Qb) \mapsto \int_{\Tbb^d} e^{i\innerprod{\kb}{\thetab}} (\Psi^{-1}+Q)^{-1} \d\nu_\Nb
	\end{equation}
	exist and are continuous in $\Dscr$. The partials evaluated at a point are viewed as linear operators between two underlying vector spaces, and continuity is understood with respect to norms.
	
	
	To ease the notation, for $(\Psi,\Qb)\in \Dscr$ let $\Phi(\Psi,\Qb):=(\Psi^{-1}+Q)^{-1}$. Consider the partial derivative w.r.t. the first argument
	\begin{align}
	D_1 \fb_{\kb}(\Psi&,\Qb):\nonumber \\& \delta\Psi \mapsto \int_{\Tbb^d} e^{i\innerprod{\kb}{\thetab}} \Phi(\Psi,\Qb) \, \Psi^{-1} \, \delta\Psi \, \Psi^{-1} \, \Phi(\Psi,\Qb) \, \d\nu_\Nb.\nonumber
	\end{align}
	Let a sequence $\{(\Psi_j,\Qb_j)\} \subset \Dscr$ converge in the product topology to $(\Psi,\Qb)\in \Dscr$, that is, $\Psi_j\to\Psi$ and $\Qb_j\to\Qb$ in respective norms. We need to show that
	\[D_1 \fb_{\kb}(\Psi_j,\Qb_j) \to D_1 \fb_{\kb}(\Psi,\Qb)\]
	in the operator norm. Indeed, we have
	\begin{equation}
	\begin{split}
	& \| D_1 \fb_{\kb}(\Psi_j,\Qb_j) - D_1 \fb_\kb(\Psi,\Qb) \| \\
	:= & \sup_{\|\delta\Psi\|=1} \| D_1 \fb_{\kb}(\Psi_j,\Qb_j;\delta\Psi)-D_1 \fb_\kb(\Psi,\Qb;\delta\Psi) \| \\
	= & \sup_{\|\delta\Psi\|=1} \| \int_{\Tbb^d} e^{i\innerprod{\kb}{\theta}}  [\, \Phi(\Psi_j,\Qb_j) \, \Psi_j^{-1} \, \delta\Psi \, \Psi_j^{-1} \, \Phi(\Psi_j,\Qb_j)  \\ & \hspace{0.4cm}  - \Phi(\Psi,\Qb) \, \Psi^{-1} \, \delta\Psi \, \Psi^{-1} \, \Phi(\Psi,\Qb) \, ] \d\nu_\Nb  \| \\
	\leq & \sup_{\|\delta\Psi\|=1}  \| \, \Phi(\Psi_j,\Qb_j) \, \Psi_j^{-1} \, \delta\Psi \, \Psi_j^{-1} \, \Phi(\Psi_j,\Qb_j)\\ & \hspace{0.4cm} - \Phi(\Psi,\Qb) \, \Psi^{-1} \, \delta\Psi \, \Psi^{-1} \, \Phi(\Psi,\Qb) \, \| \to 0.
	\end{split}
	\end{equation}
	The limit tends to zero because
	\begin{align}
	&  \| \, \Phi(\Psi_j,\Qb_j) \, \Psi_j^{-1} \, \delta\Psi \, \Psi_j^{-1} \, \Phi(\Psi_j,\Qb_j) \nonumber \\ &\hspace{0.4cm} - \Phi(\Psi,\Qb) \, \Psi^{-1} \, \delta\Psi \, \Psi^{-1} \, \Phi(\Psi,\Qb)\, \| \nonumber \\
	\leq &  \| \, \Phi(\Psi_j,\Qb_j) \, \Psi_j^{-1} \, \delta\Psi \, \Psi_j^{-1} \, \Phi(\Psi_j,\Qb_j)- \Phi(\Psi_j,\Qb_j) \, \Psi_j^{-1} \nonumber \\ &   \hspace{0.4cm}\times  \delta\Psi \, \Psi^{-1} \, \Phi(\Psi,\Qb)  +\Phi(\Psi_j,\Qb_j) \, \Psi_j^{-1} \, \delta\Psi \, \Psi^{-1} \, \Phi(\Psi,\Qb)\nonumber \\ & \hspace{0.4cm} - \Phi(\Psi,\Qb) \, \Psi^{-1} \, \delta\Psi \, \Psi^{-1} \, \Phi(\Psi,\Qb)\, \|\nonumber  \\
	\leq & \, \| \Phi(\Psi_j,\Qb_j) \, \Psi_j^{-1} \| \, \|\delta\Psi\| \, \underbrace{\|\Psi_j^{-1} \, \Phi(\Psi_j,\Qb_j) - \Psi^{-1} \, \Phi(\Psi,\Qb)\|}_{\to\, 0} \nonumber \\
	& + \underbrace{\| \Phi(\Psi_j,\Qb_j) \, \Psi_j^{-1} - \Phi(\Psi,\Qb) \, \Psi^{-1} \|}_{\to\, 0} \, \|\delta\Psi\| \, \| \Psi^{-1} \, \Phi(\Psi,\Qb) \| ,
	\end{align}
	the quantity $\|\Phi(\Psi_j,\Qb_j)\Psi_j^{-1}\|$ is bounded on $\Tbb_\Nb^d$ due to Lemma \ref{lem_bounded_below}, and we are taking the supremum over $\|\delta\Psi\|=1$.
	
	For the partial derivative of $\fb$ w.r.t. the second argument, we have
	\begin{equation}\label{f_2partial}
	D_2 \fb_\kb(\Psi,\Qb):\,\delta\Qb \mapsto -\int_{\Tbb^d} e^{i\innerprod{\kb}{\thetab}} \Phi(\Psi,\Qb) \, \delta Q \, \Phi(\Psi,\Qb) \, \d\nu_\Nb.
	\end{equation} 
	One can show the continuity of the second partial in a similar way to that for $D_1 \fb_\kb$.
	
	The same argument can be extended in a trivial manner to prove the continuity of higher-order derivatives, because the expression of $\fb_\kb$ in \eqref{f_kb} involves only rational operations on its arguments $(\Psi,\Qb)$.
\end{pf}

\end{document}